\documentclass[11pt]{amsart}
\usepackage{amssymb,amsmath,amsthm}
\usepackage{graphicx}
\usepackage{booktabs,longtable}
\usepackage[margin=1in]{geometry}
\usepackage[hidelinks]{hyperref}
\hypersetup{pdftitle={Does every non-trivial knot have segment number three?},pdfauthor={Makoto Ozawa}}

\newtheorem{theorem}{Theorem}[section]
\newtheorem{proposition}[theorem]{Proposition}
\newtheorem{lemma}[theorem]{Lemma}
\newtheorem{corollary}[theorem]{Corollary}
\newtheorem{conjecture}[theorem]{Conjecture}
\theoremstyle{definition}
\newtheorem{definition}[theorem]{Definition}
\newtheorem{remark}[theorem]{Remark}
\newtheorem{question}[theorem]{Question}

\newtheorem{example}[theorem]{Example}

\newcommand{\col}{\mathrm{col}}

\title{Does every non-trivial knot have segment number three?}
\author{Makoto Ozawa}
\address{Department of Natural Sciences, Faculty of Arts and Sciences, Komazawa University, 1-23-1 Komazawa, Setagaya-ku, Tokyo, 154-8525, Japan}
\email{w3c@komazawa-u.ac.jp}
\subjclass[2020]{Primary 57K10; Secondary 57K14}
\keywords{knot, diagram, segment number, bridge number, signature, Fox colouring}

\begin{document}
\begin{abstract}
The segment number $s(K)$ of a knot $K$ is the minimal number of arcs into which a diagram of $K$ can be cut so that no arc crosses itself and, for every pair of arcs, one of them passes over the other at all of their common crossings; every non-trivial knot has $s(K)\ge3$.
It was conjectured that the bridge number satisfies $b(K)\le s(K)$, and the torus knot $5_1$ was proposed as a candidate for $s(K)=4$.
We show that every non-trivial two-bridge knot has segment number three: a signed double spiral of half-circles around a straight arc realises every even continued fraction.
We also show that for every $n\ge2$ there are infinitely many knots with $s(K)=3$ and bridge number $n$; for $n\ge3$ they are $3$-segment plat presentations obtained by repeatedly appending one braid, and their bridge number is detected by Fox $3$-colourings.
A computer census gives $3$-segment presentations for every prime knot with at most eight crossings and for $210$ of the $249$ prime knots with at most ten crossings. We report the smallest presentations found and distinguish verified minima from results that depend on numerical knot identification.
We conjecture that every non-trivial knot has segment number three.
\end{abstract}
\maketitle

\section{Introduction}

Let $K$ be a knot in $S^3$ and $D$ a diagram of $K$ on $S^2$.
Following \cite{Oz}, an \emph{$n$-partitioned diagram} is a diagram $D$ together with $n$ points of $D$, distinct from the crossings, which cut $D$ into $n$ arcs $e_1,\dots,e_n$, called \emph{segments}.

\begin{definition}[\cite{Oz}]\label{def:segment}
An $n$-partitioned diagram $D=e_1\cup\dots\cup e_n$ is an \emph{$n$-segment presentation} of $K$ if
\begin{enumerate}
\item no segment $e_i$ has a self-crossing, and
\item for every pair $i\ne j$, either $e_i$ passes over $e_j$ at every crossing of $e_i$ and $e_j$, or $e_i$ passes under $e_j$ at every such crossing (possibly there are none).
\end{enumerate}
The \emph{segment number} $s(K)$ is the minimal $n$ such that $K$ admits an $n$-segment presentation.
\end{definition}

It is shown in \cite[Proposition 3.1]{Oz} that $K$ is non-trivial if and only if $s(K)\ge 3$, and that in a $3$-segment presentation of a non-trivial knot the three segments are cyclically stacked: $e_1$ is over $e_2$, $e_2$ is over $e_3$ and $e_3$ is over $e_1$ \cite[Proposition 3.2]{Oz}.
This ``three-way deadlock'' is the motivating picture of \cite{Oz}, where the following were posed.

\begin{conjecture}[{\cite[Conjecture 4.2]{Oz}}]\label{conj:old}
For every knot $K$, $b(K)\le s(K)$, where $b(K)$ denotes the bridge number.
\end{conjecture}

\begin{question}[{\cite[Problem 4.3]{Oz}}]\label{prob:old}
Identify a knot $K$ with $s(K)=4$; the torus knot $5_1$ was suggested as a candidate.
\end{question}

(The trivial knot has segment number one.) Our main result answers the question of the title for the two-bridge knots, including the twist knots and the torus knots $T(2,q)$. In particular, $s(5_1)=3$.

\begin{theorem}\label{thm:twobridge-intro}
Every non-trivial two-bridge knot $K$ has $s(K)=3$.
More precisely, if $K$ has Schubert fraction $p/q=[2\varepsilon_1n_1,\dots,2\varepsilon_{2g}n_{2g}]$ with $n_i\ge1$ and $\varepsilon_i=\pm1$, then $K$ has a $3$-segment presentation with
$\sum_{i}n_i+2\sum_{1\le r\le s\le g}n_{2r-1}n_{2s}$
crossings.
\end{theorem}

Every non-trivial two-bridge knot has such an even continued fraction expansion (Section~\ref{sec:51}).
The presentation is a \emph{signed double spiral}: two spirals of half-circles wind around a straight arc (Figure~\ref{fig:signed-spiral-planar}); viewed from a different direction, each crossing with the straight arc becomes one full twist, and the diagram becomes Conway's normal form (Figures~\ref{fig:signed-double-spiral-braid} and~\ref{fig:signed-spiral-local}).
Since the unknotting number, the four-ball genus and half the absolute value of the signature of $T(2,q)$ are equal to $(q-1)/2$, these invariants are unbounded on knots with $s(K)=3$ and give no lower bound for $s(K)$ beyond $3$.

\begin{theorem}\label{thm:bridge}
For every $n\ge2$ there are infinitely many knots $K$ with $s(K)=3$ and $b(K)=n$.
In particular, Conjecture~\ref{conj:old} is false, and the bridge number is not bounded by any function of the segment number.
\end{theorem}

For $n=2$ this follows from Theorem~\ref{thm:twobridge-intro}. For $n\ge3$ the knots are plat closures of $2n$-braids built from a $16$-crossing $3$-segment presentation which is a $6$-plat (Figure~\ref{fig:K40}). Two operations preserve the $3$-segment property: appending one $18$-crossing braid $G$ raises the bridge number by one, so that $n-3$ copies of $G$ give bridge number $n$; and inserting a power $L^k$, $k$ even, of a $6$-crossing braid into the initial $6$-plat keeps the bridge number fixed and changes the knot, which gives infinitely many knots for each $n$. The lower bound for the bridge number comes from the elementary inequality $\col_3(K)\le 3^{b(K)}$ for the number of Fox $3$-colourings. The smallest knot with $s(K)=3$ and $b(K)\ge4$ that we know has a $3$-segment presentation with $14$ crossings (Figure~\ref{fig:K14}).

These results suggest the opposite of the expectation in \cite[Section~4.4]{Oz}.

\begin{conjecture}[Universality of the deadlock]\label{conj:main}
Every non-trivial knot $K$ satisfies $s(K)=3$.
\end{conjecture}

By Proposition~\ref{prop:curves} below, Conjecture~\ref{conj:main} is a statement about plane curves: \emph{every knot type is represented, up to mirror image, by a generic closed curve which can be cut by three points into three simple arcs, the crossings being resolved by the cyclic rule}; one of the three arcs may even be taken to be a straight line segment (Proposition~\ref{prop:normal}).
For a non-trivial knot $K$ let $c_3(K)$ be the minimal number of crossings of a $3$-segment presentation of $K$, with $c_3(K)=\infty$ if $s(K)>3$.

\begin{theorem}[Census]\label{thm:census}
Every prime knot with at most eight crossings has segment number $3$.
Of the $249$ prime knots with at most ten crossings, the $210$ knots of Table~\ref{tab:census} have segment number $3$; they comprise all $95$ two-bridge knots and $115$ of the $154$ knots of bridge number three.
Table~\ref{tab:census} gives a presentation and hence an upper bound for $c_3(K)$ for each of these knots. Theorem~\ref{thm:enum} gives $c_3(8_{15})=27$ and $29\le c_3(8_{19})\le30$.
\end{theorem}

The $35$ presentations of the prime knots with at most eight crossings are drawn in Figures~\ref{fig:gallery} and~\ref{fig:gallery-large}.
The census is obtained by random searches and by an exhaustive enumeration of $3$-segment presentations with up to $26$ crossings (Section~\ref{sec:evidence}); all presentations are given explicitly and certified independently of the searches (Section~\ref{sec:method}).

\subsection*{What is known}
Besides the classes of Theorems~\ref{thm:twobridge-intro}, \ref{thm:bridge} and~\ref{thm:census}, $s(K)=3$ holds for the weaving knots $W(3,q)$ with $q\ge2$ and $3\nmid q$, whose standard alternating diagrams are themselves $3$-segment presentations, so that $c_3(K)=c(K)$ (Proposition~\ref{prop:braidlike}), and for the composite knots $3_1\#3_1$ and $3_1\#3_1^*$ (Example~\ref{ex:square}). The prime knots with at most ten crossings for which $s(K)$ is unknown are $9_{16}$ and $38$ knots with ten crossings, all of bridge number three (Section~\ref{sec:evidence}); among them is the torus knot $10_{124}=T(3,5)$, and for torus knots $T(p,q)$ with $3\le p<q$ the segment number is known only for $8_{19}=T(3,4)$.

On knots with $s(K)=3$, the bridge number (Theorem~\ref{thm:bridge}) and the unknotting number and the four-ball genus (Theorem~\ref{thm:twobridge-intro}, already for $T(2,q)$) are unbounded, and every pair of integer values of $a_2$ and $v_3/6$, the Vassiliev invariants of degrees $2$ and $3$, is realised by a two-bridge knot (Proposition~\ref{prop:a2v3-realize}). On the other hand, $c_3(K)$ can be much larger than $c(K)$: $c_3(K)\ge 3u(K)$ by Corollary~\ref{cor:c3}.

The paper is organised as follows.
Section~\ref{sec:prelim} contains structural facts on $3$-segment presentations.
Theorem~\ref{thm:twobridge-intro} is proved in Section~\ref{sec:51} and Theorem~\ref{thm:bridge} in Section~\ref{sec:bridge}.
Section~\ref{sec:twist} gives constructions beyond two-bridge knots, and Section~\ref{sec:examples} the examples mentioned above.
Section~\ref{sec:invariants} expresses the Conway coefficient $a_2$ and the Vassiliev invariant $v_3$ as signed counts of crossings of a $3$-segment presentation.
Section~\ref{sec:evidence} describes the census and the enumeration, Section~\ref{sec:method} the certification of the data, and Section~\ref{sec:questions} lists open questions.

\section{Structure of \texorpdfstring{$3$}{3}-segment presentations}\label{sec:prelim}

Diagrams are polygonal.
We use the convention of \cite{Oz} that a $3$-segment presentation $D=e_1\cup e_2\cup e_3$ is traversed in the order $e_1,e_2,e_3$, with partition points $p_1=e_1\cap e_2$, $p_2=e_2\cap e_3$, $p_3=e_3\cap e_1$; thus $e_1$ runs from $p_3$ to $p_1$, $e_2$ from $p_1$ to $p_2$ and $e_3$ from $p_2$ to $p_3$.

\subsection{Tri-simple curves}

\begin{definition}
A \emph{tri-simple curve} is a generic closed polygonal curve $C\subset\mathbb{R}^2$ together with three points $p_1,p_2,p_3\in C$, distinct from the double points, which cut $C$ into three arcs $e_1,e_2,e_3$ (in this cyclic order), each of which is simple.
\end{definition}

Every double point of a tri-simple curve lies on two different arcs, so the following \emph{cyclic rule} is well defined: at a double point of $e_i$ and $e_{i+1}$ (indices mod $3$), $e_i$ is the over-strand.
Denote the resulting knot by $K_+(C)$, and the knot obtained by the opposite rule by $K_-(C)$.

\begin{proposition}\label{prop:curves}
\begin{enumerate}
\item $K_-(C)$ is the mirror image of $K_+(C)$.
\item A knot $K$ satisfies $s(K)\le 3$ if and only if $K$ or its mirror image is $K_+(C)$ for some tri-simple curve $C$.
\item $s(K)=s(K^*)=s(-K)$, where $K^*$ denotes the mirror image and $-K$ the reverse.
\end{enumerate}
\end{proposition}

\begin{proof}
(1) The two rules differ at every crossing, and changing every crossing of a diagram yields a diagram of the mirror image.
(2) The underlying curve of a $3$-segment presentation is tri-simple by condition (1) of Definition~\ref{def:segment}.
If the presentation represents a non-trivial knot, its digraph is an oriented $3$-cycle \cite[Proposition 3.2]{Oz}, so its crossings follow one of the two cyclic rules.
The trivial knot is $K_+(C)$ for a round circle $C$.
Conversely, $K_\pm(C)$ with the partition $e_1,e_2,e_3$ satisfies Definition~\ref{def:segment}.
(3) Changing all crossings of a $3$-segment presentation gives a $3$-segment presentation of $K^*$, and the definition does not involve orientations.
\end{proof}

Thus $s(K)\le 3$ is a property of a \emph{projection} together with three cut points: the crossing information is forced up to mirror image.

\begin{proposition}[Normal form]\label{prop:normal}
If $s(K)\le 3$, then $K$ or $K^*$ has a $3$-segment presentation in which $e_3$ is the straight line segment from $p_2=(0,0)$ to $p_3=(1,0)$.
The same holds with $e_3$ replaced by $e_1$ or $e_2$.
\end{proposition}

\begin{proof}
For a simple polygonal arc $e$ in $\mathbb{R}^2$ there is an orientation-preserving PL homeomorphism of $\mathbb{R}^2$ carrying $e$ onto a straight segment (PL Schoenflies theorem applied to a regular neighbourhood of $e$).
Its product with the identity of the vertical coordinate is an orientation-preserving homeomorphism of $\mathbb{R}^3$ which respects the over/under information of the diagram, hence preserves the knot type; after a small perturbation the image is again a generic polygon.
\end{proof}

Hence, up to mirror image, the knots with $s(K)\le 3$ are exactly the knots obtained from a pair of simple arcs, $e_1$ from $(1,0)$ to a point $p_1$ and $e_2$ from $p_1$ to $(0,0)$, by adding the \emph{bar} $e_3=[0,1]\times\{0\}$, which passes over $e_1$ and under $e_2$, while $e_1$ passes over $e_2$.

\subsection{Elementary bounds}
Recall that $c_3(K)$ denotes the minimal number of crossings of a $3$-segment presentation of a non-trivial knot $K$.
For a $3$-segment presentation $D=e_1\cup e_2\cup e_3$ we write $C_{ij}$ for the set of crossings between $e_i$ and $e_j$ and $c_{ij}=|C_{ij}|$.
We say that $D$ \emph{follows the cyclic rule} if $e_1$ passes over $e_2$, $e_2$ over $e_3$ and $e_3$ over $e_1$ at all their crossings. By \cite[Proposition 3.2]{Oz}, every $3$-segment presentation of a non-trivial knot follows the cyclic rule or its reverse, and by reversing the traversal or taking a mirror image we may assume the former; but a $3$-segment presentation in the sense of Definition~\ref{def:segment} may also have totally ordered segments, and then it represents the trivial knot.

\begin{proposition}\label{prop:unknotting}
Let $D$ be a $3$-segment presentation of $K$ which follows the cyclic rule.
Changing the crossings of $C_{ij}$, for any one pair $\{i,j\}$, yields a descending diagram of the trivial knot.
Hence
\[
u(K)\le \min\{c_{12},\,c_{23},\,c_{13}\},
\]
and the same bound holds for the ascending number of \cite{Oz-asc}.
\end{proposition}

\begin{proof}
After changing the crossings of $C_{13}$, the segment $e_1$ passes over $e_3$, so $e_1$ is over $e_2$ and $e_3$, and $e_2$ is over $e_3$.
Since no segment has a self-crossing, the resulting diagram, traversed from $p_3$ along $e_1,e_2,e_3$, meets every crossing first as an over-crossing; that is, it is descending, and a descending diagram represents the trivial knot.
Thus $K$ is unknotted by the $c_{13}$ crossing changes at $C_{13}$, and $u(K)\le c_{13}$.
Changing $C_{12}$ gives the order $e_2>e_3>e_1$, descending from $p_1$, and changing $C_{23}$ gives $e_3>e_1>e_2$, descending from $p_2$; hence $u(K)\le c_{12}$ and $u(K)\le c_{23}$.
The diagrams obtained are descending, so the same changes also show that the ascending number of $K$ is at most $\min\{c_{12},c_{23},c_{13}\}$.
\end{proof}

\begin{corollary}\label{cor:c3}
For every non-trivial knot $K$,
\[
c_3(K)\ \ge\ \max\{\,c(K),\ 3u(K)\,\}.
\]
More precisely, every $3$-segment presentation of $K$ satisfies $c_{ij}\ge u(K)$ for all three pairs $\{i,j\}$.
\end{corollary}

\begin{proof}
The second statement is Proposition~\ref{prop:unknotting}.
Summing over the three pairs, a $3$-segment presentation has $c_{12}+c_{23}+c_{13}\ge 3u(K)$ crossings, and it has at least $c(K)$ crossings since it is a diagram of $K$.
\end{proof}

For example $c_3(5_1)\ge 6$, $c_3(7_1)\ge 9$ and $c_3(8_{19})\ge 9$, since $u(5_1)=2$ and $u(7_1)=u(8_{19})=3$.

A \emph{bigon} of a diagram is a face whose boundary contains exactly two crossings.
A \emph{twist region} is a maximal chain of bigons stacked end to end, or a single crossing not incident to any bigon; $t(D)$ denotes the number of twist regions of $D$.

\begin{proposition}\label{prop:bigon}
Let $K$ be a non-trivial knot with $s(K)=3$, and let $D$ be a $3$-segment presentation of $K$ with $c_3(K)$ crossings.
Then:
\begin{enumerate}
\item every bigon of $D$ has a partition point on its boundary;
\item $D$ has at most $6$ bigons;
\item $t(D)\ge c(D)-6$; in particular every twist region of $D$ has at most $7$ crossings.
\end{enumerate}
\end{proposition}

\begin{proof}
(1) Let $F$ be a bigon with crossings $x$ and $y$.
If neither side of $F$ contains a partition point, each side lies in a single segment, say $e_i$ and $e_j$; here $i\ne j$ because no segment has a self-crossing.
By condition (2) of Definition~\ref{def:segment}, $e_i$ passes over $e_j$ at both $x$ and $y$ (or under at both), so a Reidemeister II move across $F$ removes $x$ and $y$.
The resulting diagram is again a $3$-segment presentation with the same partition points, because removing crossings cannot violate the conditions of Definition~\ref{def:segment}.
It has $c_3(K)-2$ crossings, a contradiction.

(2) A partition point lies in the interior of an edge of the diagram, and an edge lies on the boundary of at most two faces.
By (1) each bigon contains one of the three partition points on its boundary, so there are at most $3\cdot 2=6$ bigons.

(3) A twist region with $k$ crossings contains at least $k-1$ bigons, and distinct twist regions share no bigon.
Summing over the twist regions $R$ of $D$, with $k_R$ crossings and $b_R\ge k_R-1$ bigons,
\[
c(D)=\sum_R k_R\le \sum_R (b_R+1)\le 6+t(D),
\]
by (2). A twist region with $k$ crossings has $k-1\le6$ bigons, so $k\le7$.
\end{proof}

The moves used in the proof of Proposition~\ref{prop:bigon} are the only local moves available within $3$-segment presentations.

\begin{lemma}\label{lem:rigid}
Let $D$ be a $3$-segment presentation which follows the cyclic rule and in which no edge of the diagram contains two partition points; the latter holds whenever $D$ represents a non-trivial knot.
Neither a Reidemeister move III across a triangular face of $D$ which keeps the partition points on their edges, nor a slide of a partition point across a crossing, transforms $D$ into a $3$-segment presentation with the same cyclic rule.
\end{lemma}

\begin{proof}
Label each passage of the knot through a crossing by the index of the segment containing it.
Since $D$ follows the cyclic rule, at every crossing the labels $o$ of the over-passage and $u$ of the under-passage satisfy $u\equiv o+1 \pmod 3$; neither move changes which strand is the over-strand at a crossing.
A slide of a partition point across a crossing $x$ changes the label of exactly one of the two passages at $x$, by $\pm1$, so afterwards $u\not\equiv o+1$ at $x$.

Now let $D'$ be obtained by a move III across a triangular face with crossings $x,y,z$, whose sides lie on strands $A$ (from $x$ to $y$), $B$ (from $x$ to $z$) and $C$ (from $y$ to $z$); along each strand the two crossings are passed in the opposite order after the move.
By assumption no edge of $D$ contains two partition points. (If an edge contains two partition points, one segment has no crossings, and $D$ is a diagram in which one of the other two segments passes over the other at every crossing, a diagram of the trivial knot; so the assumption holds when $D$ represents a non-trivial knot.)
If the side on $A$ contains no partition point, the labels of $A$ at $x$ and $y$ do not change.
If it contains one, the labels of $A$ at $x$ and $y$ are $s$ and $s+1$ in the order of $A$, and after the move the crossing passed first has label $s$; so the labels at $x$ and $y$ change by $\delta_A(x)=-\delta_A(y)=\pm1$.
Define $\delta_B$ and $\delta_C$ in the same way.
If $D'$ is a $3$-segment presentation, the condition $u\equiv o+1$ at each crossing gives $\delta_A(x)=\delta_B(x)$, $\delta_A(y)=\delta_C(y)$ and $\delta_B(z)=\delta_C(z)$, hence
\[
\delta_A(x)=\delta_B(x)=-\delta_B(z)=-\delta_C(z)=\delta_C(y)=\delta_A(y)=-\delta_A(x),
\]
so all $\delta$ vanish and no side of the face contains a partition point.
Then each side lies in a single segment; the three segments are pairwise distinct since the sides cross each other, so they are $e_1,e_2,e_3$, and by the cyclic rule none of the three strands passes over both of the others.
A Reidemeister move III requires such a strand, a contradiction.
\end{proof}

Consequently, the only Reidemeister moves between $3$-segment presentations with the same cyclic rule which keep the partition points on their edges are moves I and II. Hence a presentation in which no move I or II removes crossings cannot be simplified by a sequence of such moves which never increases the number of crossings and keeps the partition points and the cyclic rule.

Thus a minimal $3$-segment presentation is ``almost twist-free'': almost every crossing forms a twist region of its own.
For example, for odd $n\ge5$ no minimal diagram of $T(2,n)$ is a $3$-segment presentation, since $3u(T(2,n))=3(n-1)/2>n=c(T(2,n))$ (Corollary~\ref{cor:c3}).

\section{Two-bridge knots}\label{sec:51}
\label{subsec:all-two-bridge}

We prove Theorem~\ref{thm:twobridge-intro}. For a sequence of non-zero integers we write $[a_1,\dots,a_n]=a_1+1/(a_2+1/(\cdots+1/a_n))$, and $S(p,q)$ for the two-bridge knot with Schubert fraction $p/q$ \cite[Chapter~12]{BZ}.

\begin{theorem}\label{thm:all-two-bridge}
Let $a_i=2\varepsilon_i n_i$, where $g\ge1$, $n_i\ge1$, $\varepsilon_i=\pm1$, and
$1\le i\le 2g$. Let $p/q=[a_1,\ldots,a_{2g}]$, with the plus sign at each
step of the continued fraction. The signed double spiral described below
is a three-segment presentation of the two-bridge knot of slope $p/q$,
up to mirror image. Its number of crossings is
\[
  \sum_{i=1}^{2g} n_i
    +2\!\!\sum_{1\le r\le s\le g} n_{2r-1}n_{2s}.
\]
\end{theorem}

\paragraph{The planar double spiral.}

Put $e_1=[0,B]$ on the horizontal axis, oriented to the right, and
write $p_3=0$, $p_1=B$. Place marked points in $(0,B)$ in the order
$x^{n_1}w^{n_2}\cdots x^{n_{2g-1}}w^{n_{2g}}$. Each marked point in the
$i$th block receives the direction $\varepsilon_i$, with $+1$ meaning
upward in the original projection. Denote by $d(t)$ the direction
assigned to a marked point $t$.

The arc $e_2$ starts at $B$ and visits the $x$ points in decreasing order.
Just after each visit there is a turn at an exterior point to the right
of $B$. Use strictly increasing exterior coordinates, then close the arc
with one semicircle to $p_2<0$. The arc $e_3$ starts at $p_2$ and visits
the $w$ points in decreasing order. Just after each visit it turns at
an exterior point in $(p_2,0)$, with strictly increasing exterior
coordinates, and ends at $0$. Each consecutive pair of vertices is
joined by a semicircle. For either arc, take the incoming semicircle at
its first marked point on side $-d(t)$ and alternate sides throughout.
If consecutive marked points of that arc require opposite directions,
insert one additional exterior vertex between them. The resulting
extra semicircle has both ends outside $[0,B]$. The endpoint $0$ may
serve as the last exterior turn of $e_3$.
Otherwise, join the last exterior vertex to $0$ by a short semicircle.

\begin{figure}[htbp]
\centering
\includegraphics[width=.8\linewidth]{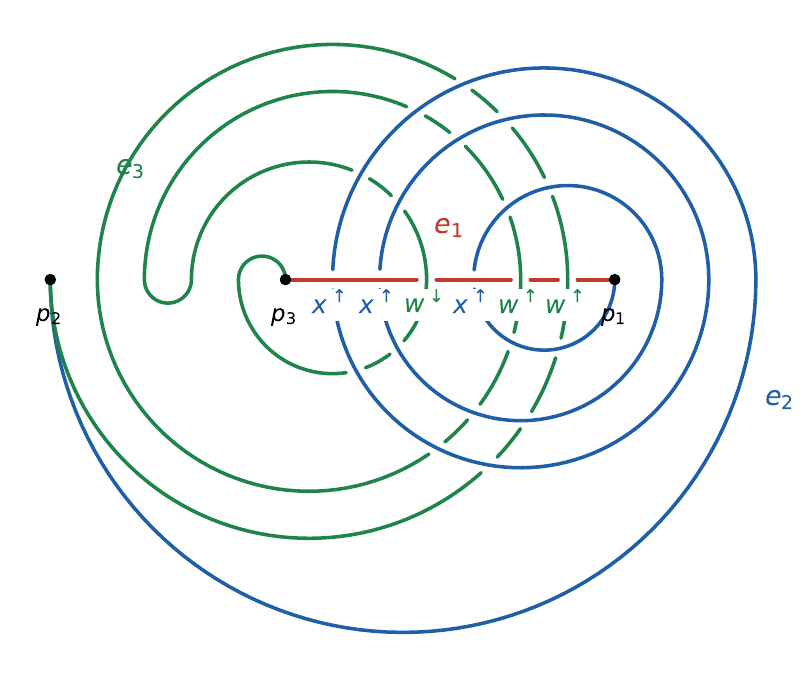}
\caption{The signed double spiral for $[a_1,a_2,a_3,a_4]=[4,-2,2,4]$, a $3$-segment presentation of $10_{16}$ with $22$ crossings. The straight arc $e_1$ (red) runs from $p_3$ to $p_1$; the marked points carry the word $x^2w\,x\,w^2$, and the arrows give the directions $\varepsilon_i$ of the four blocks. The arc $e_2$ (blue) spirals out of $p_1$ and visits the $x$ points from right to left; the arc $e_3$ (green) spirals into $p_3$ and visits the $w$ points from right to left. The small hairpin of $e_3$ on the left is the extra exterior vertex which reverses the direction between its two blocks. At every crossing $e_1$ is over $e_2$, $e_2$ is over $e_3$ and $e_3$ is over $e_1$.}\label{fig:signed-spiral-planar}
\end{figure}

The diameters of any two long semicircles of the same arc are nested
or share an end. The short exterior semicircles have disjoint interiors
on the axis and do not interleave any diameters of their own arc.
The closing diameter of $e_2$ contains all the other diameters.
Thus there are no self-crossings. Perturb corners at the three partition
points if necessary. Every $x$ and every $w$ contributes one crossing
with $e_1$. For a fixed pair $x,w$, the long semicircles based at the
two points interleave exactly when $x<w$. When they interleave, there
is one crossing in each half-plane. The short semicircles contribute
none: those on the right are outside all diameters of $e_3$, and those
on the left do not interleave the closing diameter of $e_2$.
This gives the asserted count and a generic three-segment diagram (Figure~\ref{fig:signed-spiral-planar}).

\paragraph{Straightening the passes.}

Use coordinates $(x,y,z)$ in space, with the spiral drawn in the
$xy$-plane and $z$ recording over and under. Reserve disjoint small
neighborhoods of the marked points of $e_1$. In each such neighborhood
the only pieces of the spiral are the axis and the short transverse
piece that crosses it. Tilt each transverse piece slightly so that it
travels from larger to smaller $x$. It still crosses $y=0$ in its
assigned direction.

Lift $e_1$ to height zero. Away from the neighborhoods and the partition
points, put $e_2$ at height $+2$ and $e_3$ at height $-2$. At an $x$ point,
dip the short piece of $e_2$ to height $-1$ as it crosses $y=0$.
At a $w$ point, raise the short piece of $e_3$ to height $+1$ as it
crosses $y=0$. Across each short piece, prescribe the two crossings of
height zero to lie on opposite sides of $y=0$. The height can be
chosen piecewise linear and smoothed. Connect the heights near the
partition points in small neighborhoods disjoint from the other
projected arcs. These choices realize the original cyclic crossings.

We use the following elementary replacement. Let a pass $\alpha$ of the knot lie at height $H$ over a simple planar path $\gamma$, and let $\gamma'$ be another simple planar path with the same end points, agreeing with $\gamma$ on short end collars. Suppose that, away from those collars, every point of the knot other than $\alpha$ which lies over $\gamma\cup\gamma'$ has height less than $H$. Then $\alpha$ can be replaced by the arc at height $H$ over $\gamma'$, by an isotopy fixed near the end collars. Indeed, the vertical region above height $H$ over $\gamma\cup\gamma'$ meets no other part of the knot: lift $\alpha$, except for short end collars, vertically above all remaining strands; at this height the two simple paths are isotopic in the plane relative to their end points; and lower the new pass vertically along $\gamma'$. The end collars are kept in small neighbourhoods of the end points throughout. The same argument applies with ``below'' in place of ``above'', and heights greater than $H$.

Apply this operation to the parts of $e_2$ and $e_3$ between the short
transverse pieces. Here is how to choose the paths so the operation
can be done successively.

\begin{itemize}
\item For two consecutive visits of one arc, use a path moving
monotonically in $x$ in a thin ribbon about the axis between the
two reserved neighborhoods. Use the same choice from $p_1$ to the
first $x$, from the last $x$ towards $p_2$, and from the last $w$ to
$p_3$. When the path to $p_2$ leaves $[0,B]$, follow a small
neighborhood of the axis until it reaches the left endpoint.
\item Outside the reserved neighborhoods, the original projected
pieces of that same arc are disjoint from the portion of the axis
in the interior of the interval in question. In the end pieces,
the part of the axis between $p_2$ and the last $x$ is also free of
other visits of $e_2$. By compactness, the ribbon can be chosen thin
enough to avoid all those pieces. If the interval includes marks of
the other arc, those pieces are on the correct side in height: below
$+2$ or above $-2$, respectively. Near an endpoint, route the new path
through the side of the reserved neighborhood facing into its
interval. The short crossing piece is monotone in $x$ there and
the next pass leaves from the other end of that short piece.
Previously replaced passes of the same arc lie in disjoint
interval interiors. Thus, outside the end collars, every other point of the knot over the old and the new path of a pass of $e_2$ has height less than $+2$, and over those of a pass of $e_3$ height greater than $-2$: the parts of the other two arcs lie at heights in $[-2,+1]$, respectively $[-1,+2]$, and the parts of the same arc at height $\pm2$ are avoided in projection. So the hypotheses of the replacement hold at each step.
\item Replace the initial pass of $e_3$ last. Run it from the rightmost
$w$ to the right of $B$, then around the ribbon to $p_2$ on its exterior,
reading that pass backwards. The new path can be disjoint in projection
from every other part of $e_3$, which is already inside the ribbon
apart from its short end collars. Put the long return part at height
$-3$, keeping the transitions to the height $-2$ in the two end
regions. In the new view, this is an outside closing arc from the
bottom right strand to the top left strand.
\end{itemize}

\begin{figure}[htbp]
\centering
\includegraphics[width=.95\linewidth]{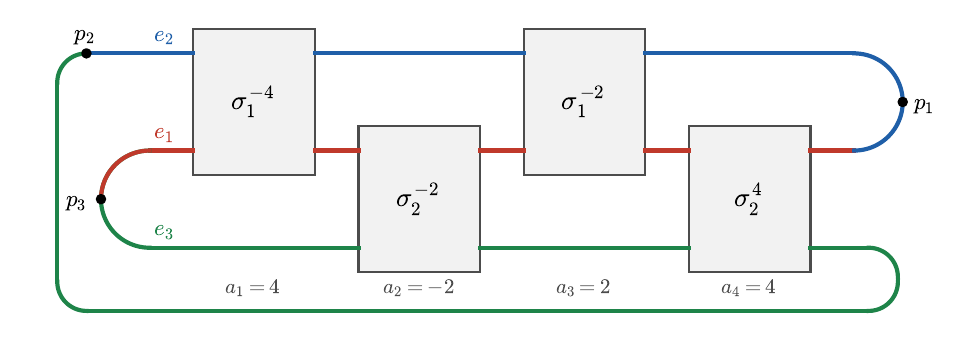}
\caption{The diagram of Figure~\ref{fig:signed-spiral-planar} after straightening the passes, projected to the $(x,z)$-plane. From top to bottom the strands are $e_2,e_1,e_3$; each box is a full-twist region, $\sigma_1^{-a_i}$ for an $x$-block and $\sigma_2^{a_i}$ for a $w$-block, read from left to right. The partition point $p_3$ caps $e_1$ and $e_3$ on the left, $p_1$ caps $e_2$ and $e_1$ on the right, and the outside arc through $p_2$ joins the remaining two ends: the even-length Conway closure.}\label{fig:signed-double-spiral-braid}
\end{figure}

View the result in the $(x,z)$-plane from $y<0$. Inside the braid box,
from top to bottom the strands are normally $e_2,e_1,e_3$.
At the left end $p_3$ caps the middle and bottom strands.
At the right end $p_1$ caps the top and middle strands. These two
caps can be drawn outside the braid box in their small end
neighborhoods. The remaining ends are joined by the outside arc
through $p_2$. This is the even-length Conway closure (Figure~\ref{fig:signed-double-spiral-braid}).

\paragraph{The braid word.}

\begin{figure}[htbp]
\centering
\includegraphics[width=.95\linewidth]{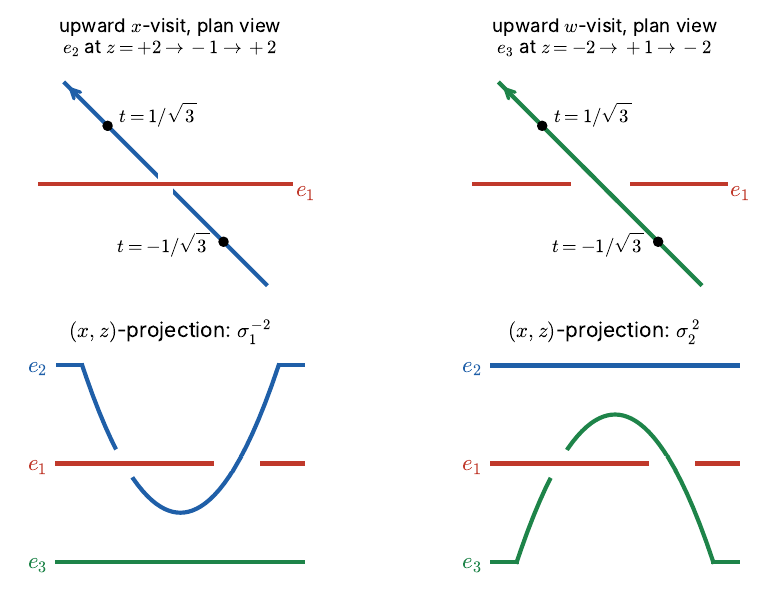}
\caption{The local model at an upward $x$-visit (left) and an upward $w$-visit (right). Top: the plan view, in which the transverse piece meets the plane $z=0$ at $t=\pm1/\sqrt3$, once on each side of $e_1$. Bottom: the projection to the $(x,z)$-plane seen from $y<0$, where the nearer strand is over; the piece crosses $e_1$ twice with the same sign, giving the full twists $\sigma_1^{-2}$ and $\sigma_2^{2}$.}\label{fig:signed-spiral-local}
\end{figure}

Read the braid from left to right, with $\sigma_i$ denoting a crossing in
which the strand in position $i$ passes over the strand in position $i+1$.
The positions are numbered from top to bottom just before the crossing.
Choose $\delta,\eta>0$ small. In a neighborhood of an upward visit at
$x=x_0$ we may take the oriented transverse piece to be
\[
 (x,y,z)=
 \begin{cases}
 (x_0-\delta t,\ \eta t,\ 2-3\phi(t)), &\text{at an $x$-visit (on $e_2$)},\\
 (x_0-\delta t,\ \eta t,\ -2+3\phi(t)), &\text{at a $w$-visit (on $e_3$)},
 \end{cases}
 \qquad
 \phi(t)=1-t^2,\quad -1\le t\le1 .
\]
Here $t$ increases along the oriented piece. In the braid projection we
encounter $t=1/\sqrt3$ at the left crossing, with $y>0$ (farther from
the viewer), and $t=-1/\sqrt3$ at the right crossing, with $y<0$
(nearer to the viewer). The axis $e_1$ has $y=0$ throughout. Thus the
two crossings are as follows (Figure~\ref{fig:signed-spiral-local}).

\begin{center}
\begin{tabular}{lcc}
 & Left crossing ($y>0$) & Right crossing ($y<0$)\\[2pt]
 $x$-visit & $e_2$ is above and passes under $e_1$
           & $e_1$ is above and passes under $e_2$\\
 $w$-visit & $e_1$ is above and passes over $e_3$
           & $e_3$ is above and passes over $e_1$
\end{tabular}
\end{center}
Here \emph{above} refers to the position in the $(x,z)$-projection just
before the crossing, and \emph{over} refers to the smaller $y$-coordinate.
The two strands exchange positions at the first crossing. Hence the
$x$-visit gives $\sigma_1^{-2}$ and the $w$-visit gives $\sigma_2^{2}$.
Each is one full twist, consisting of two half-twists of the same sign.
For a downward visit replace $\eta t$ by $-\eta t$; this reverses both
crossings. The neighborhoods have disjoint $x$-intervals, so changing a
visit direction changes only its own full twist.
Therefore the
word is
\[
  \sigma_1^{-a_1}\,\sigma_2^{a_2}\,\sigma_1^{-a_3}
  \cdots\sigma_2^{a_{2g}} .
\]
After reflecting the diagram, the word becomes
$\sigma_1^{a_1}\sigma_2^{-a_2}\sigma_1^{a_3}\cdots\sigma_2^{-a_{2g}}$,
with the same closure. This is Conway's normal form for
$[a_1,\ldots,a_{2g}]$ in the indicated convention
\cite[Chapter~12, Sections~A--B]{BZ}. In slope coordinates its successive
blocks act by the shears $r\mapsto r+a$ and $r\mapsto r/(br+1)$,
respectively. Thus its slope is the image of $\infty$ under the product
of the matrices
\[
 \begin{pmatrix}1&a_1\\0&1\end{pmatrix}
 \begin{pmatrix}1&0\\a_2&1\end{pmatrix}
 \cdots
 \begin{pmatrix}1&0\\a_{2g}&1\end{pmatrix}
 =
 \prod_{i=1}^{2g}\begin{pmatrix}a_i&1\\1&0\end{pmatrix},
\]
which evaluates to $[a_1,\ldots,a_{2g}]$. Reflecting the original
projection across $y=0$ with its crossing
heights fixed reverses every visit direction and gives the other
mirror, with the same cyclic rule. This proves the theorem.

\begin{proof}[Proof of Theorem~\ref{thm:twobridge-intro}]
For a non-trivial two-bridge knot choose a Schubert fraction $p/q$ with
$p>1$ odd and $q$ even, $0<|q|<p$. One obtains such a representative by
adding or subtracting $p$ from an odd denominator in $(0,p)$.
Apply the Euclidean algorithm with each quotient chosen as the
nearest even integer. At each step $(P,Q)$ is a coprime pair of
opposite parity with $|P|>|Q|$. An odd integer cannot be $P/Q$,
since that would force $|Q|=1$ and $P$ odd. Hence there is no tie,
the quotient is nonzero and even, and $|P-aQ|<|Q|$.
The parities interchange at each step. The final denominator is
$\pm1$, so the last step is even-numbered. This gives a nonzero
even expansion of even length. The theorem gives $s(K)\le3$. The lower bound $s(K)\ge3$ for non-trivial knots is \cite[Proposition~3.1]{Oz}.
\end{proof}

\subsection{Torus knots \texorpdfstring{$T(2,q)$}{T(2,q)} and twist knots}\label{sec:T2q}
The crossing count of Theorem~\ref{thm:all-two-bridge} depends on the expansion. For $[2a,2b]$ with $a,b>0$ it is $2ab+a+b$; for the twist knot $C(n,2)=S(2n+1,2)$ (up to mirror image), with $c=n+2$ crossings and expansion $[2a,\pm2]$, $a=\lceil n/2\rceil$, it is $3a+1\le 3c/2-1/2$.
The torus knot $T(2,2k+1)$ has the expansion $[2,-2,\dots,2,-2]$ of length $2k$, and the signed double spiral is then a presentation with $k^2+3k$ crossings in which the two spirals are interleaved without hairpins.
The trefoil is the case $k=1$; for $k=2$ we obtain a presentation of $5_1$ with $10$ crossings, so $5_1$ is not an example for Question~\ref{prob:old}.
One crossing can be saved.

\begin{proposition}[computer-assisted]\label{prop:T2q}
For $1\le k\le 15$ the signed Gauss words $D_k$ given by an explicit formula in the supplementary data (\texttt{family.py}) are $3$-segment presentations of $T(2,2k+1)$ with $(c_{12},c_{23},c_{13})=(k,k,k^2+k-1)$, hence with $k^2+3k-1$ crossings.
In particular, $c_3(3_1)=3$, $c_3(5_1)\le9$, $c_3(7_1)\le17$ and $c_3(9_1)\le27$.
\end{proposition}

\begin{proof}
For each $k\le15$ a script of the supplementary data checks that $D_k$ is realised by a diagram in the sphere (by the Euler characteristic), that the three segments are simple and follow the cyclic rule, and that the fundamental group of the complement, computed by SnapPy \cite{SnapPy}, has the presentation $\langle a,b\mid a^2=b^{2k+1}\rangle$; torus knots are determined by their groups up to mirror image \cite{GL,Whitten}. The bounds follow from the crossing counts, and $c_3(3_1)=3$ since $c(3_1)=3$.
\end{proof}

The words $D_k$ describe $e_1$ and $e_3$ as two interleaved families of nested half-circles, as in the double spiral. The targeted searches for small $k$ are discussed in Section~\ref{sec:evidence}.
By Corollary~\ref{cor:c3} below, $c_3(T(2,2k+1))\ge 3k$. Thus $c_3$ grows at most linearly on twist knots, and at most quadratically on the knots $T(2,q)$; we do not know whether $c_3(T(2,q))/q$ is bounded (Question~\ref{q:T2q}).

\section{Segment number does not bound bridge number}\label{sec:bridge}

For an odd prime $p$, a \emph{Fox $p$-colouring} of a diagram assigns elements of $\mathbb{Z}/p$ to its arcs so that $2x_{\mathrm{over}}=x_{\mathrm{in}}+x_{\mathrm{out}}$ at every crossing.
The number $\col_p(K)$ of $p$-colourings is a knot invariant, equal to $p^{1+r}$ with $r=\dim_{\mathbb{Z}/p}H_1(\Sigma_2(K);\mathbb{Z}/p)$, where $\Sigma_2(K)$ is the double branched cover \cite{Przytycki}.

\begin{lemma}\label{lem:col-bridge}
For every knot $K$ and odd prime $p$, $\col_p(K)\le p^{b(K)}$.
\end{lemma}

\begin{proof}
Let $b=b(K)$ and take a diagram coming from a $b$-bridge position: the knot is the union of $b$ over-bridges and $b$ under-bridges, the over-strand at every crossing belongs to an over-bridge, and the over-bridges do not cross each other.
Each over-bridge lies in a single arc of the diagram.
Travelling along an under-bridge from its initial point, which lies on the arc of the preceding over-bridge, the colour of each arc is determined by the colour of the previous arc and the colour of the over-strand.
Hence a colouring is determined by the colours of the $b$ arcs containing the over-bridges.
\end{proof}

\begin{figure}[ht]
\centering
\includegraphics[width=0.5\textwidth]{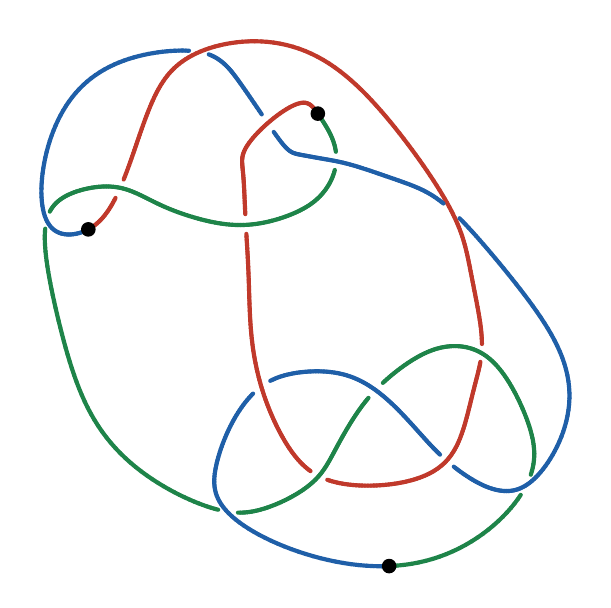}
\caption{A $3$-segment presentation with $14$ crossings of the knot $K=\texttt{K14n14349}$, which has $\col_3(K)=3^4$, hence $b(K)\ge 4$. The segments $e_1$ (red), $e_2$ (blue), $e_3$ (green) are separated by the three black dots, and $e_1>e_2>e_3>e_1$. The drawing has the signed Gauss code of the polygon of Table~\ref{tab:coordsK14}, and is isotopic to it on $S^2$.}
\label{fig:K14}
\end{figure}

\begin{proposition}\label{prop:K14}
The knot $K$ of Figure~\ref{fig:K14} has a $3$-segment presentation with $14$ crossings and $b(K)\ge4$.
\end{proposition}

\begin{proof}
Let $D$ be the polygon with the vertices of Table~\ref{tab:coordsK14} (Figure~\ref{fig:K14}).
Exact arithmetic shows that $D$ is a generic $3$-segment presentation with $14$ crossings ($c_{12}=5$, $c_{23}=5$, $c_{13}=4$) and hence $14$ arcs.
The $14\times14$ colouring system over $\mathbb{Z}/3$ has rank $10$, so $\col_3(K)=3^4$.
By Lemma~\ref{lem:col-bridge}, $b(K)\ge 4$.
On the other hand $K$ is non-trivial, since the trivial knot has only $3$ colourings; therefore $s(K)=3$.
\end{proof}

\begin{remark}
SnapPy \cite{SnapPy} identifies $K$ as the hyperbolic knot \texttt{K14n14349}, with $H_1(\Sigma_2(K))\cong\mathbb{Z}/3\oplus\mathbb{Z}/3\oplus\mathbb{Z}/21$.
This is not needed for Proposition~\ref{prop:K14}. Since $c(K)=14$, the presentation of Figure~\ref{fig:K14} has the minimal possible number of crossings; it was obtained from a $24$-crossing presentation found by the search of Section~\ref{sec:evidence} by moves that preserve the $3$-segment property and decrease the number of crossings (e.g.\ Reidemeister II moves, which are always available on a bigon between two segments).
Several other knots with $s(K)=3$ and $\col_3(K)=3^4$ occurred in our experiments, e.g.\ \texttt{K14n20076}, and a knot with $H_1(\Sigma_2)\cong\mathbb{Z}/3\oplus\mathbb{Z}/15\oplus\mathbb{Z}/15$.
\end{remark}

\subsection{Knots of every bridge number}\label{sec:every-bridge}
For a braid $\beta\in B_{2n}$ the \emph{plat closure} of $\beta$ joins the end points $2j-1$ and $2j$, $j=1,\dots,n$, by caps at both ends of $\beta$; it has a $2n$-plat diagram, so its bridge number is at most $n$. We read braids from left to right and let $\sigma_i$ be the crossing of the strands at positions $i$ and $i+1$ in which the strand moving from position $i$ to position $i+1$ passes over; $\beta^{(j)}$ denotes $\beta$ with all indices shifted by $j$. Let
\begin{align*}
L&=(\sigma_4\sigma_5^{-1})^3,\\
\beta_0^{(k)}&=\sigma_2^{-1}\sigma_3^{-1}\sigma_4^{-1}\,L^k\,\sigma_5\sigma_1^{-1}\sigma_2^{-1}\sigma_4^{-1}\sigma_3\sigma_2^{-1}\sigma_5\sigma_3\sigma_2^{-1}\sigma_3\sigma_4^{-1}\sigma_3\sigma_4^{-1}\in B_6,\\
G&=\sigma_1\sigma_2^{-1}\sigma_1\sigma_2^{-1}\sigma_1\;\sigma_3\sigma_4^{-1}\sigma_3\sigma_4^{-1}\sigma_3\;\sigma_2^{-1}\sigma_3\sigma_2^{-1}\sigma_3\;\sigma_5\sigma_4^{-1}\sigma_3\sigma_4^{-1}\in B_6,
\end{align*}
and, for $n\ge3$ and $k\ge0$,
\[
\beta_{n,k}=\beta_0^{(k)}\,G^{(2)}G^{(4)}\cdots G^{(2n-6)}\in B_{2n},
\]
where $\beta_0^{(k)}\in B_6\subset B_{2n}$. Let $K_{n,k}$ be the plat closure of $\beta_{n,k}$. It has $16+6k+18(n-3)$ crossings; $K_{3,0}$ is the knot $12a_{503}$.

\begin{figure}[ht]
\centering
\includegraphics[width=\textwidth]{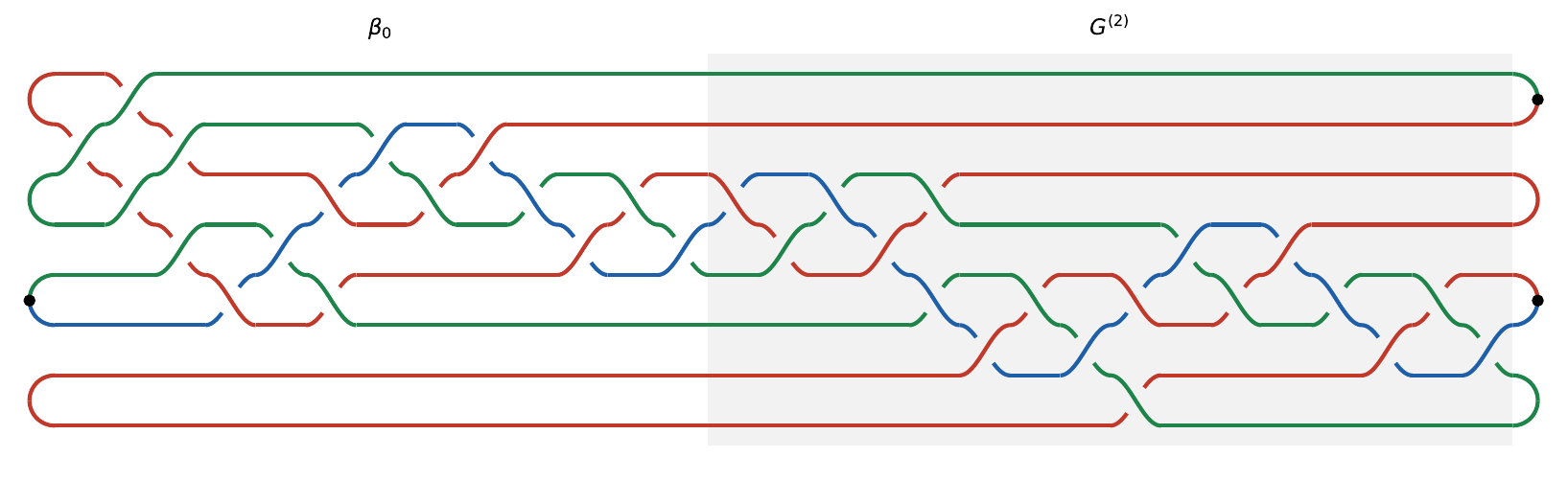}
\caption{The $8$-plat closure of $\beta_{4,0}=\beta_0^{(0)}G^{(2)}$, a $3$-segment presentation with $34$ crossings of a knot $K_{4,0}$ with $b(K_{4,0})=4$. The braid is read from left to right; the part $G^{(2)}$ is shaded. The strands are labelled $A$ (red), $B$ (blue), $C$ (green), and $A>B>C>A$ at every crossing; the three partition points (black dots) lie on the caps where the label changes.}
\label{fig:K40}
\end{figure}

\begin{theorem}\label{thm:every-bridge}
Let $n\ge3$ and let $k\ge0$ be even.
\begin{enumerate}
\item The plat closure of $\beta_{n,k}$ is a $3$-segment presentation of the knot $K_{n,k}$.
\item $\col_3(K_{n,k})=3^n$; hence $b(K_{n,k})=n$ and $s(K_{n,k})=3$.
\item $\det K_{n,k}=\bigl|(1,0)\,A^{n-3}v_k\bigr|$, where
\[
A=\begin{pmatrix}2955&1800\\1896&1155\end{pmatrix},\quad
v_0=\begin{pmatrix}-135\\-87\end{pmatrix},\quad
v_1=\begin{pmatrix}21\\13\end{pmatrix},\quad
v_{k+2}=18v_{k+1}-v_k .
\]
In particular, for fixed $n$ the determinant is strictly increasing along even $k$, and the knots $K_{n,k}$, $k$ even, are pairwise distinct.
\end{enumerate}
\end{theorem}

\begin{proof}
(1) Label the left end points of $\beta_{n,k}$ by $A,A,C,C,C,B$ followed by $A,A,\dots,A$ (positions $7,\dots,2n$), and carry the labels along the strands. A direct check shows that every crossing of $\beta_0^{(k)}$ is a crossing of two strands with different labels $X,Y$ at which $X$ passes over $Y$ for $(X,Y)\in\{(A,B),(B,C),(C,A)\}$; this holds for every $k$, since $L$ acts on three strands with the three labels at positions $4,5,6$ and is a pure braid. The right labels of $\beta_0^{(k)}$ are $C,A,A,B,C,C$. Similarly $G$, applied to strands labelled $A,B,C,C,A,A$, follows the same rule at all its $18$ crossings and ends with the labels $A,A,A,B,C,C$; so the last four labels are again $B$ preceded by $A$ and followed by $C,C$, and $G^{(2)},\dots,G^{(2n-6)}$ can be appended inductively. The right labels of $\beta_{n,k}$ are $C,A,\dots,A,B,C,C$. Thus the label changes only on the left cap $(5,6)$ and on the right caps $(1,2)$ and $(2n-3,2n-2)$; placing the three partition points there cuts the closure into three arcs labelled $A$, $B$, $C$, which follow the cyclic rule.
The closure is a knot: this is checked directly for $n=3$, and the strands of $G$ end in the order $(1,5,6,2,3,4)$, so that in passing from $\beta_{n-1,k}$ to $\beta_{n,k}$ the right cap joining the local positions $1,2$ is replaced by an arc through the new left cap $(5,6)$, and the number of components does not change.

(2) Consider the colours of the strands of $G$ over $\mathbb{Z}$. If the left colours are $(a,a,b,b,c,c)$, the differences of the colours joined by the three right caps of $G$ are
\[
105a-120b+15c,\qquad -825a+960b-135c,\qquad 720a-840b+120c,
\]
all divisible by $3$. Hence every $3$-colouring of $K_{n-1,k}$, which has equal colours $a,a,b,b$ at the right end points $2n-5,\dots,2n-2$ of $\beta_{n-1,k}$, extends to a $3$-colouring of $K_{n,k}$ for each of the three values of the colour $c$ of the new left cap, and distinct pairs give distinct colourings. So $\col_3(K_{n,k})\ge3\col_3(K_{n-1,k})$. The braid $L$ acts on the $3$-colourings of three strands with order $4$, so $\col_3(K_{3,k})$ depends only on $k\bmod4$; a computation gives $\col_3(K_{3,0})=\col_3(K_{3,2})=3^3$. Thus $\col_3(K_{n,k})\ge3^n$, and by Lemma~\ref{lem:col-bridge} and the $2n$-plat diagram, $3^n\le\col_3(K_{n,k})\le3^{b(K_{n,k})}\le3^n$. So $b(K_{n,k})=n$, and $s(K_{n,k})=3$ by~(1), since $K_{n,k}$ is non-trivial.

(3) Normalise the colours over $\mathbb{Z}$ so that the colours of the closed caps are equal. Then the colours of the last four strands of $\beta_{n-1,k}$ can be written as $(r+t,t,0,r)$ up to a common translation. In the system of cap equations, with the constant colourings removed, the coefficient of the colour $c$ of the new left cap in the equation of the first right cap of $G$ is $15$ (see~(2)); we use this equation, with pivot $15$, to eliminate $c$. Computing the new pair $(r',t')$ gives $(r',t')^{\mathsf T}=\frac1{15}A\,(r,t)^{\mathsf T}$, and the pivot $15$ cancels this denominator in the determinant of the colouring matrix, which gives the induction step $v\mapsto Av$. For $\beta_0^{(k)}$ the same computation gives the vectors $v_k$. The braid $L$ acts on the Fox colourings of its three strands by the integer matrix
\[
R=\begin{pmatrix}9&-4&-4\\-4&1&4\\-12&4&9\end{pmatrix},
\]
and $\ell R=\ell$ for $\ell=(1,-1,1)$, so that $\ker\ell$ is invariant under $R$; the restriction of $R$ to $\ker\ell$ has characteristic polynomial $\lambda^2-18\lambda+1$. If the left caps of $\beta_0^{(k)}$ carry the colours $(0,0,a,a,b,b)$, the colours entering $L$ at the positions $4,5,6$ are $b\,w$ with $w=(1,2,1)^{\mathsf T}$, which do not involve $a$; so, after the elimination with the first right cap, $v_k=\Lambda R^k w$ for a fixed linear map $\Lambda\colon\mathbb{Q}^3\to\mathbb{Q}^2$. Since $w\in\ker\ell$, the Cayley--Hamilton theorem for $R|_{\ker\ell}$ gives $(R^2-18R+I)R^kw=0$, and applying $\Lambda$ yields $v_{k+2}=18v_{k+1}-v_k$ for all $k\ge0$. This gives the formula, which we have also checked against the determinant of the colouring matrix for $3\le n\le 6$ and $0\le k\le6$. The entries of $A$ are positive, $v_k$ has positive, strictly increasing entries for $k\ge1$, and $v_2=(513,321)^{\mathsf T}$ exceeds $|v_0|$ entrywise; hence $\det K_{n,k}$ is strictly increasing along even $k$.
\end{proof}

For $k=0$ the determinants are $135$, $555{,}525$, $2{,}283{,}177{,}375$ and $9{,}383{,}734{,}018{,}125$ for $n=3,4,5,6$, with $16$, $34$, $52$ and $70$ crossings. All presentations with $n\le6$ and $k\in\{0,2\}$ were also checked independently: the plat closures were built from the braid words, the $3$-segment property was confirmed by a search over all triples of partition points, and $\col_3$ and $\det$ were computed from the colouring matrices (supplementary data, \texttt{bridge\_family/}).

\begin{proof}[Proof of Theorem~\ref{thm:bridge}]
For $n=2$ this is Theorem~\ref{thm:twobridge-intro}; for instance the $4$-plat closures of $(\sigma_2\sigma_3^{-1})^{3k}\sigma_2$, $k\ge1$, are $3$-segment presentations of two-bridge knots with determinant $F_{6k+2}$. For $n\ge3$ it is Theorem~\ref{thm:every-bridge}.
\end{proof}

\section{Beyond two-bridge knots}\label{sec:twist}

We give two constructions of $3$-segment presentations which are not restricted to two-bridge knots: twisting along a circle around two partition points, and braid-like diagrams.

\subsection{A height model}
Let $D=e_1\cup e_2\cup e_3$ be a $3$-segment presentation which follows the cyclic rule.
A \emph{height model} of $D$ is a knot in $\mathbb{R}^3$ projecting onto $D$ in which
$e_1$ lies in the plane $z=0$, $e_3$ lies in the plane $z=1$, $e_2$ is a graph over its projection lying below $z=0$ at the crossings with $e_1$ and above $z=1$ at the crossings with $e_3$, and the three partition points are replaced by vertical segments.
Since $e_3$ is above $e_1$, $e_1$ above $e_2$ and $e_2$ above $e_3$ at every crossing, a height model represents the knot $K(D)$ given by $D$.
The height of $e_2$ away from its crossings is at our disposal.

\subsection{The spiral lemma}
Let $0<r_0<r_1$, let $B=\{r<r_0\}$ be a disc and $A=\{r_0\le r\le r_1\}$ an annulus around it, in polar coordinates $(r,\theta)$ of the plane, and let $\tau$ be the Dehn twist of the plane supported in $A$, $\tau(r,\theta)=(r,\theta+2\pi\varphi(r))$, where $\varphi=1$ on $[0,r_0]$ and $\varphi=0$ on $[r_1,\infty)$.

\begin{lemma}[Spiral lemma]\label{lem:spiral}
Let $D$ be a $3$-segment presentation following the cyclic rule. Suppose that no partition point lies in $A$, every component of $e_3\cap A$ is an arc joining the two boundary circles of $A$, and no crossing between $e_1$ and $e_2$ lies in $A$.
For $k\in\mathbb{Z}$ let $D^{(k)}$ be obtained from $D$ by replacing $e_3\cap A$ with $\tau^k(e_3\cap A)$ and keeping $e_1$ and $e_2$.
Then $D^{(k)}$ is a $3$-segment presentation.
Moreover, choose a height model of $D$ in which $e_2\cap A$ lies above $z=1+2\varepsilon$, and let $c_\pm=\partial B\times\{1\pm\varepsilon\}$.
Then $K(D^{(k)})$ is obtained from $K(D)$ by $-k$ full twists along $c_-$ and $k$ full twists along $c_+$, with the sign convention fixed in the proof.
\end{lemma}

\begin{proof}
The arc $\tau^k(e_3)$ is simple, since $\tau$ is a homeomorphism of the plane and $e_3\cap A$ is mapped into $A$.
The only new crossings lie in $A$, between the spiral and the arcs of $e_1\cap A$ and $e_2\cap A$; the cyclic rule makes $e_3$ pass over $e_1$ and under $e_2$ there, so $D^{(k)}$ is a $3$-segment presentation.
The height model of $D$ can be chosen with $e_2\cap A$ above $z=1+2\varepsilon$ because $e_2$ has no crossing with $e_1$ in $A$.
We also choose it so that, over $B$, the knot meets the slab $1-\varepsilon\le z\le1+\varepsilon$ only in $e_3\cap B$, at height $1$, and in finitely many vertical segments contained in pairwise disjoint thin vertical cylinders over points of $B$ away from the crossings. This is possible: $e_2$ changes its height only on short subarcs away from the crossings and the partition points, and we may make $e_2$ cross the slab vertically on such subarcs; the vertical segments replacing the partition points are vertical by definition. Thus over $B$ the part of the model in each of the two slabs $1-\varepsilon\le z\le1$ and $1\le z\le1+\varepsilon$ is a trivial tangle of vertical segments.
Let $\psi\colon\mathbb{R}\to[0,1]$ be supported in $(1-\varepsilon,1+\varepsilon)$ with $\psi(1)=1$, increasing below $1$ and decreasing above $1$, and let
\[
\Phi_t(r,\theta,z)=\bigl(r,\ \theta+2\pi t k\,\varphi(r)\psi(z),\ z\bigr),\qquad 0\le t\le 1 .
\]
This is an ambient isotopy of $\mathbb{R}^3$.
In $A\times\mathbb{R}$ the height model meets the slab $1-\varepsilon<z<1+\varepsilon$ only in $e_3$, at height $1$, so $\Phi_1$ replaces $e_3\cap A$ by $\tau^k(e_3\cap A)$ and fixes the rest of the model over $A$ and outside $B\cup A$.
Over $B$, $\Phi_1$ rotates the plane $z=1$ by $2\pi k$, so it fixes $e_3\cap B$; it acts only on the vertical segments of the model which cross the slab over $B$, and on $B\times[1-\varepsilon,1]$ and on $B\times[1,1+\varepsilon]$ it is $k$, respectively $-k$, full twists of these segments.
Here a \emph{positive full twist} of the vertical segments over $B$ is the rotation by $2\pi$ in the direction of increasing $\theta$ as $z$ increases, and a full twist along $c_\pm$ is the corresponding twist of the knot along the disc $B\times\{1\pm\varepsilon\}$; with this convention $\Phi_1$ adds $k$ positive full twists in $B\times[1-\varepsilon,1]$ and $-k$ in $B\times[1,1+\varepsilon]$, since $\psi$ increases from $0$ to $1$ on the first slab and decreases on the second. Undoing these twists of trivial tangles, that is, applying $-k$ full twists along $c_-$ and $k$ full twists along $c_+$, gives a height model of $D^{(k)}$. We use the same convention for the twists along $c_-$ in Proposition~\ref{prop:corner}.
\end{proof}

\begin{remark}\label{rem:lk}
The linking numbers of $c_\pm$ with $K(D)$ lie in $\{-1,0,1\}$.
Indeed, cut the height model along the cylinder $\partial B\times\mathbb{R}$; the ends of the pieces over $B$ lie on $e_1$ (below $c_-$), on $e_3$ (between $c_-$ and $c_+$) or on $e_2$ (above $c_+$), so a piece crosses the disc bounded by $c_-$ algebraically a non-zero number of times only if it contains $p_3$ (from $e_3$ to $e_1$), or $p_1$ (from $e_1$ to $e_2$), and each partition point occurs once; the same argument applies to $c_+$, with $p_1$ and $p_2$.
In particular the spiral lemma never twists along a circle of linking number $2$, such as the braid axis of a closed $2$-braid, and it does not directly give the torus knots $T(2,q)$; these are obtained by the double spiral of Section~\ref{sec:51}.
\end{remark}

\subsection{Corner twists}
The spiral lemma applies to a $3$-segment presentation following the cyclic rule when two partition points lie on a common face.
Let $D$ be a $3$-segment presentation, and suppose that two partition points lie on the closure of a common face $F$ of $D$ (a component of the complement of $D$ in the plane).
Since the indices are cyclic, we may call them $p_a$ and $p_b$, where $e_m$ runs from $p_a$ to $p_b$, $e_{m-1}$ ends at $p_a$ and $e_{m+1}$ starts at $p_b$.
Let $\alpha$ be an arc in $F\cup\{p_a,p_b\}$ joining $p_a$ to $p_b$, and let $N$ be a thin regular neighbourhood of $\alpha$, so that $D\cap N$ consists of four arcs without crossings, two of $e_m$ and one each of $e_{m-1}$ and $e_{m+1}$, each joining a partition point to $\partial N$.
The \emph{corner circle} of $(D,\alpha)$ is the unknot $c$ which projects onto $\partial N$, passes over the two arcs of $e_m$ and under the arcs of $e_{m-1}$ and $e_{m+1}$.
It bounds a disc meeting $K(D)$ in two points of opposite signs, so $\operatorname{lk}(c,K(D))=0$.

\begin{proposition}\label{prop:corner}
Let $D$ be a $3$-segment presentation following the cyclic rule, with $N$ crossings, and let $c$ be a corner circle of $(D,\alpha)$.
For every $k\in\mathbb{Z}$, the knot obtained from $K(D)$ by $k$ full twists along $c$ has a $3$-segment presentation with $N+3|k|$ crossings.
\end{proposition}

\begin{proof}
First slide $p_b$ along $\alpha$ to a point $q$ close to $p_a$: replace $e_m$ by $e_m$ followed by a parallel copy of the sub-arc of $\alpha$ from $p_b$ to $q$, and $e_{m+1}$ by a second parallel copy, on the other side, followed by $e_{m+1}$.
Since $\alpha$ lies in a face, no crossing is created, the new partition point is $q$, and the new presentation $D'$ represents the same knot (the two copies form a finger without crossings, which can be retracted in a height model); the corner circle is carried along.
Now let $B$ be a small disc containing $p_a$ and $q$ and $A$ a thin annulus around it, both inside the image of $N$, so that $D'\cap(B\cup A)$ consists of four arcs without crossings, each crossing $A$ once: the ends of $e_{m-1}$ and of $e_m$ at $p_a$, the end of $e_m$ at $q$ and the beginning of $e_{m+1}$ at $q$.
Relabel the segments cyclically so that $e_{m-1},e_m,e_{m+1}$ become $e_3,e_1,e_2$; a cyclic relabelling preserves the cyclic rule, so we are in the situation of Lemma~\ref{lem:spiral}, whose hypotheses clearly hold.
For $j\in\mathbb{Z}$ the presentation $D'^{(j)}$ has $N+3|j|$ crossings, since the spiral crosses each of the other three arcs $|j|$ times.
In a height model of $D'$ in which $e_2$ lies above $z=1+2\varepsilon$ over $B\cup A$, the knot crosses the slab $1-\varepsilon<z<1+\varepsilon$ over $B$ only in the vertical segment at $p_a$ (from $e_3$ at height $1$ to $e_1$ at height $0$) and in the vertical segment at $q$ (from $e_1$ at height $0$ to $e_2$).
Thus the disc bounded by $c_+$ meets the knot once and is disjoint from $c_-$, so twisting along $c_+$ does not change the pair $(K(D'),c_-)$. By Lemma~\ref{lem:spiral} with $j=-k$, the knot $K(D'^{(-k)})$ is obtained from $K(D)$ by $k$ full twists along $c_-$, and $D'^{(-k)}$ has $N+3|k|$ crossings.
The circle $c_-=\partial B\times\{1-\varepsilon\}$ lies above the two arcs of $e_1$ (the old $e_m$) and below the arcs of $e_3$ and $e_2$ (the old $e_{m-1}$ and $e_{m+1}$), so it is the corner circle.
\end{proof}

\begin{remark}\label{rem:corner}
A search of the presentations of $38$ knots of Table~\ref{tab:census} found $41$ corner circles, in the presentations of $32$ of these knots, each giving an infinite family of distinct knots. Many of these families contain knots which are not two-bridge; for instance the corner twists of $8_{21}$ include $6_2$, $9_{42}$ and $10_{133}$. For the twist knots, the corner twists of the standard diagram of $3_1$ give presentations of $C(n,2)$ with $3+3|n|$ crossings, the link $3_1\cup c$ being the Whitehead link.
The presentations of $9_2$ and $10_{10}$ in Table~\ref{tab:census} were obtained in this way, and so were the first presentations of $10_1$ and $10_{133}$, later improved by the enumeration of Section~\ref{sec:evidence}.
By Remark~\ref{rem:lk}, corner twists, like all applications of the spiral lemma, twist along circles of linking number at most one.
\end{remark}

\subsection{Braid-like presentations}
A different mechanism gives $3$-segment presentations which are minimal diagrams.
Consider a diagram containing three arcs $A$, $B$, $C$ which contain all the crossings, such that the complement of $A\cup B\cup C$ in the knot consists of three arcs without crossings, and such that for each pair of $A,B,C$ one of them passes over the other at all their common crossings. Placing one partition point on each of the three complementary arcs gives a $3$-segment presentation (the cyclic order in which the knot traverses $A,B,C$ only decides whether the cyclic rule or its reverse holds, which is irrelevant for the segment number).
For $F_0=0$, $F_1=1$, $F_{n+1}=F_n+F_{n-1}$, the \emph{Fibonacci two-bridge knot} is $C(1,\dots,1)=S(F_{n+1},F_n)$ with $n$ entries; it is a knot if and only if $F_{n+1}$ is odd, that is, $n\not\equiv 2\pmod 3$. The \emph{weaving knot} $W(3,q)$ is the closure of the $3$-braid $(\sigma_1\sigma_2^{-1})^q$; it is a knot if and only if $3\nmid q$.

\begin{proposition}\label{prop:braidlike}
For $n\ge3$ with $n\not\equiv2\pmod3$, and for $q\ge2$ with $3\nmid q$, the standard alternating diagrams of $S(F_{n+1},F_n)$ and of $W(3,q)$ are $3$-segment presentations. Hence these knots have segment number $3$ and $c_3(K)=c(K)$, equal to $n$ and $2q$, respectively.
\end{proposition}

\begin{proof}
Let $\sigma_i$ denote the crossing in which the strand at position $i$ passes over the strand at position $i+1$, and $\sigma_i^{-1}$ the crossing in which it passes under it. We follow the strands, not the positions. Starting with the strands $A,B,C$ at positions $1,2,3$, the word $\sigma_2\sigma_1^{-1}\sigma_2\sigma_1^{-1}\sigma_2\sigma_1^{-1}$ gives:
\begin{center}
\begin{tabular}{ccccccc}
\toprule
crossing & $1$ & $2$ & $3$ & $4$ & $5$ & $6$\\
\midrule
letter & $\sigma_2$ & $\sigma_1^{-1}$ & $\sigma_2$ & $\sigma_1^{-1}$ & $\sigma_2$ & $\sigma_1^{-1}$\\
order before & $ABC$ & $ACB$ & $CAB$ & $CBA$ & $BCA$ & $BAC$\\
over-strand & $B>C$ & $C>A$ & $A>B$ & $B>C$ & $C>A$ & $A>B$\\
\bottomrule
\end{tabular}
\end{center}
After six crossings the order and the next letter are the initial ones, so every crossing of a word $\sigma_2\sigma_1^{-1}\sigma_2\cdots$ of any length satisfies $A>B>C>A$. Each of $A,B,C$ is a strand of the braid, hence has no self-crossing.
The Fibonacci knot $S(F_{n+1},F_n)$ is the four-plat closure of the first $n$ letters of this word on four strands, the fourth strand being without crossings: the top caps join the positions $(1,2)$ and $(3,4)$, and the bottom caps join $(1,2)$, $(3,4)$ if $n$ is odd and $(1,4)$, $(2,3)$ if $n$ is even. The caps together with the fourth strand form three arcs without crossings joining the ends of $A$, $B$, $C$.
Similarly, $W(3,q)$ is the closure of the first $2q$ letters of the word $\sigma_1\sigma_2^{-1}\sigma_1\sigma_2^{-1}\cdots$, to which the same table applies after the symmetry exchanging the positions $1$ and $3$; the three closing arcs have no crossings.
In both cases the remark above gives a $3$-segment presentation; the diagrams are reduced and alternating, so they have the minimal number of crossings.
\end{proof}

For example, the Fibonacci knots with $n=3,4,6,7,9,10$ are $3_1$, $4_1$, $6_3$, $7_7$, $9_{31}$ and $10_{45}$, and $W(3,2)$, $W(3,4)$, $W(3,5)$, $W(3,7)$ are $4_1$, $8_{18}$, $10_{123}$ and \texttt{K14a19470}; this agrees with the values of $c_3$ found by the enumeration of Section~\ref{sec:evidence} (Table~\ref{tab:census}).
The mechanism is rigid: two consecutive equal letters $\sigma_i^{\pm1}\sigma_i^{\pm1}$ make the same two strands cross twice with the over-strand exchanged. Accordingly, among the closed $3$-braids given by words of length at most $14$, those whose crossings are consistent along the strands in this sense and which are knots are exactly the unknot and the weaving knots $W(3,q)$; two-bridge knots with longer twist regions require other constructions, such as the double spiral of Section~\ref{sec:51}.

\section{Examples}\label{sec:examples}

\subsection{Signature and four-genus}\label{sec:4d}
In a random sample of $82{,}985$ tri-simple curves (Section~\ref{sec:evidence}) all knots had $|\sigma|\le4$, and all $10{,}558$ knots of a subsample for which $\tau$ was computed had $|\tau|\le2$. These apparent bounds are artefacts of uniform sampling, which strongly favours simple knots: a search which \emph{optimises} the signature quickly produced the following example, and the torus knots $T(2,q)$ have unbounded signature by Theorem~\ref{thm:twobridge-intro}.

\begin{figure}[ht]
\centering
\includegraphics[width=0.8\textwidth]{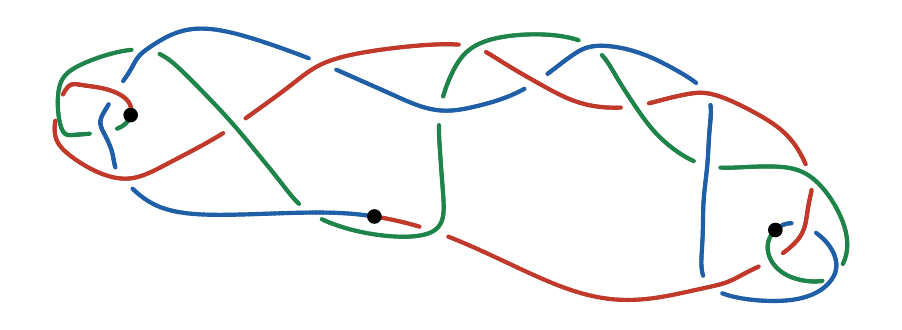}
\caption{A $3$-segment presentation with $21$ crossings of a knot with $\sigma=-6$ and $\tau=-3$. The segments $e_1$ (red), $e_2$ (blue), $e_3$ (green) are separated by the three black dots, and $e_1>e_2>e_3>e_1$. The drawing has the signed Gauss code of the polygon whose coordinates are given in the supplementary data, and is isotopic to it on $S^2$.}
\label{fig:sig6}
\end{figure}

\begin{example}\label{ex:sig}
There is a knot $K$ with $s(K)=3$, $\sigma(K)=-6$ and $\tau(K)=-3$; in particular $u(K)\ge g_4(K)\ge3$. Here $\tau$ is the Ozsv\'ath--Szab\'o concordance invariant, with $|\tau(K)|\le g_4(K)\le u(K)$ \cite{OS}.
\end{example}

\begin{proof}
Let $D$ be the diagram of Figure~\ref{fig:sig6}, with integer coordinates given in the supplementary data.
Exact arithmetic shows that $D$ is a generic $3$-segment presentation with $21$ crossings ($c_{12}=c_{23}=c_{13}=7$).
The signature of $D$, computed by the Gordon--Litherland formula \cite{GLi}, is $-6$.
The knot Floer homology, computed with Szab\'o's program, gives $\tau=-3$ and Seifert genus $7$, and shows that the knot is fibred.
Since $|\tau|\le g_4\le u$ \cite{OS}, we get $u(K)\ge g_4(K)\ge 3$.
\end{proof}

\subsection{Composite knots}\label{sec:square}
We first observe that a $3$-segment presentation can never exhibit a connected sum decomposition.

\begin{lemma}\label{lem:notvisible}
Let $D$ be a $3$-segment presentation of a knot $K$, and let $C\subset S^2$ be a circle meeting $D$ transversely in two points which are not crossings.
Then the diagram obtained from one of the two sides of $C$ (closing the arc of $D$ on that side by an arc of $C$) represents the trivial knot.
\end{lemma}

\begin{proof}
The two points of $C\cap D$ cut $K$ into two arcs $\gamma$ and $\gamma'$, and the crossings on each side of $C$ are self-crossings of the corresponding arc.
The three partition points are distributed between $\gamma$ and $\gamma'$, so one of them, say $\gamma$, contains at most one partition point.
Then $\gamma$ is contained in the union of at most two consecutive segments, and the diagram $D_\gamma$ obtained by closing $\gamma$ along $C$ is a $k$-segment presentation with $k\le 2$: its pieces are subarcs of segments of $D$, so they have no self-crossings and any two of them cross consistently.
By \cite[Proposition 3.1]{Oz}, $D_\gamma$ represents the trivial knot.
\end{proof}

Hence a $3$-segment presentation of a composite knot is never visibly composite.
Nevertheless composite knots can have segment number three.

\begin{example}\label{ex:square}
The square knot $3_1\#3_1^*$ and the granny knot $3_1\#3_1$ have segment number $3$.
\end{example}

\begin{figure}[ht]
\centering
\includegraphics[width=0.55\textwidth]{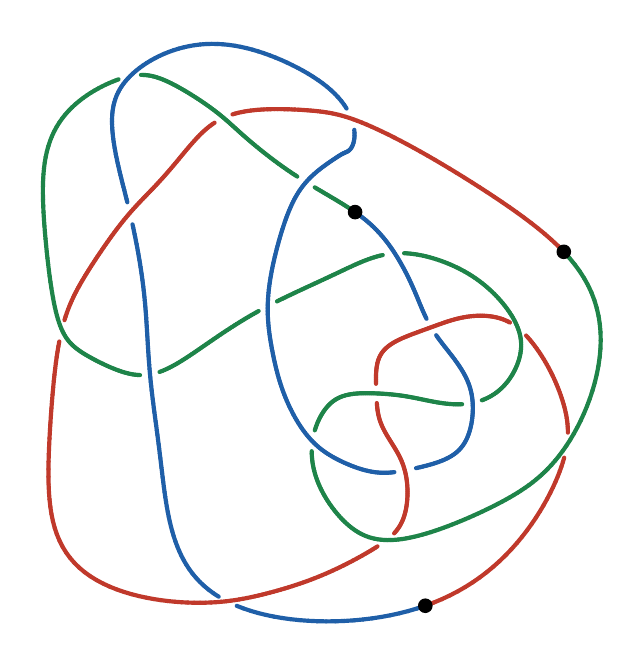}
\caption{A $3$-segment presentation with $18$ crossings of the square knot $3_1\#3_1^*$. The segments $e_1$ (red), $e_2$ (blue), $e_3$ (green) are separated by the three black dots, and $e_1>e_2>e_3>e_1$. The drawing has the signed Gauss code of the polygon of Table~\ref{tab:coordsSq}, and is isotopic to it on $S^2$.}
\label{fig:square}
\end{figure}

\begin{proof}[Proof for the square knot]
Let $D$ be the diagram of Figure~\ref{fig:square}, with vertices as in Table~\ref{tab:coordsSq}.
Exact arithmetic shows that $D$ is a generic $3$-segment presentation with $18$ crossings ($c_{12}=5$, $c_{23}=7$, $c_{13}=6$).
We exhibit an explicit sequence $D=D_0,D_1,\dots,D_{45}$ of diagrams, each given by its signed Gauss code (supplementary data), such that:
\begin{enumerate}
\item the signed Gauss code of $D_0$ is that of the integer polygon $D$;
\item each $D_{k+1}$ is obtained from $D_k$ by a single Reidemeister move performed across a face of the planar diagram: across a monogon (move~I), a bigon whose two crossings have the same over-strand (move~II), or a triangle one of whose sides passes over or under both of the others (move~III), or by the inverse of such a move;
\item the signed Gauss code of $D_{45}$ is the concatenation of two words $a\,b\,c\,a\,b\,c$ and $d\,e\,f\,d\,e\,f$ in disjoint letters, each alternating between over- and under-passages, with all crossings of the first word of one sign and all crossings of the second word of the other sign.
\end{enumerate}
A signed Gauss code determines the cyclic order of the four edges at each crossing, hence the diagram on $S^2$ up to homeomorphism, and its faces; conditions (1)--(3) are therefore finite combinatorial checks.
They were verified by a short program that uses no knot-theoretic software: it reconstructs the faces of each diagram from the code, checks that each diagram is planar (it has $n+2$ faces, where $n$ is the number of crossings), and checks that consecutive codes differ by one of the moves in (2), up to relabelling of the crossings and reversal of the orientation.
The $45$ moves consist of $4$ moves I and $9$ moves II that remove crossings, $5$ inverse moves II, and $27$ moves III; the number of crossings rises from $18$ to at most $24$ and then decreases to $6$.
By (3), $D_{45}$ is a connected sum of two standard trefoil diagrams of opposite chirality.
Hence $D$ represents $3_1\#3_1^*$.
As further checks, $D$ has signature $0$, $\col_3=3^3$, and its knot Floer homology has total rank $9$, $\tau=0$ and Seifert genus $2$, as expected for the square knot.
Since the square knot is non-trivial, $s=3$.
\end{proof}

The presentation was found by the random search of Section~\ref{sec:evidence}, as a $34$-crossing tri-simple curve, and reduced to $18$ crossings by moves preserving the $3$-segment property and the knot type.
It has no bigons.
\begin{figure}[ht]
\centering
\includegraphics[width=0.55\textwidth]{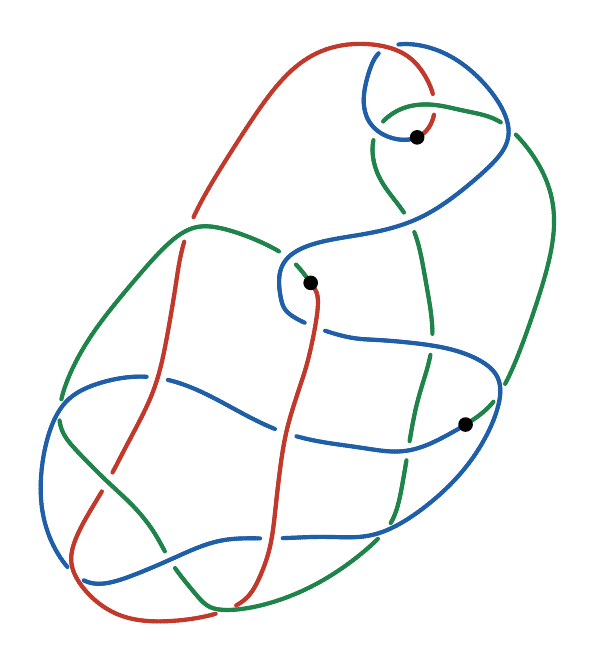}
\caption{A $3$-segment presentation with $20$ crossings of the granny knot $3_1\#3_1$. The segments $e_1$ (red), $e_2$ (blue), $e_3$ (green) are separated by the three black dots; $e_1$ passes over $e_2$ six times, $e_2$ over $e_3$ ten times, and $e_3$ over $e_1$ four times. The drawing has the signed Gauss code of the polygon given in the supplementary data, and is isotopic to it on $S^2$.}
\label{fig:granny}
\end{figure}

\begin{proof}[Proof for the granny knot]
We use the diagram $D'$ of Figure~\ref{fig:granny}, whose signed Gauss code is given in the supplementary data.
Each coloured arc is simple, and the $20$ crossings consist of $c_{12}=6$, $c_{23}=10$, $c_{13}=4$ crossings with the cyclic stacking, so $D'$ is a $3$-segment presentation; this is checked directly on the Gauss code, and the code is planar (it has $22$ faces).
We give two independent certificates that $D'$ represents $3_1\#3_1$.

\emph{(a) Re-routing of passes.}
Six times, an arc which passes entirely over (or entirely under) the rest of the diagram is removed and replaced by another such arc with the same endpoints.
In each case the remaining arc, a connected diagram on $S^2$, is unchanged: its signed Gauss code, read from the endpoint of the removed arc, is the same before and after the replacement, so its embedding is the same up to homeomorphism.
Since the removed arc can be moved freely above (or below) the rest of the diagram, the knot type is preserved.
The number of crossings goes $20\to17\to15\to12\to10\to8\to6$, and the final signed Gauss code is $(+1,-2,+3,-1,+2,-3\mid +4,-5,+6,-4,+5,-6)$, with all six crossings of the same sign: the standard connected sum of two trefoil diagrams of the same chirality.

\emph{(b) Reidemeister moves.}
As in the proof for the square knot, a sequence of $22$ Reidemeister moves performed across faces (with $10$ moves III) transforms $D'$ into the same kind of $6$-crossing diagram, and every step was verified by the program of the supplementary data.

As further checks, $D'$ has signature $\pm4$, $\col_3=3^3$, and knot Floer homology of total rank $9$ with $|\tau|=2$ and Seifert genus $2$, as expected for the granny knot.
Since $3_1\#3_1$ is non-trivial, $s=3$.
\end{proof}

The exhaustive search of Section~\ref{sec:evidence} found no $3$-segment presentation of $3_1\#3_1$ with at most $9$ crossings, which supports $c_3(3_1\#3_1)\ge10$ (determinant $9$, all cyclic classes of $(c_{12},c_{23},c_{13})$ with entries $\ge u=2$ and sum at most $9$).

\section{Invariants of \texorpdfstring{$3$}{3}-segment presentations}\label{sec:invariants}

Finite type invariants of a knot with a $3$-segment presentation can be read off from the order of the crossings along the segments. We show this for $a_2$ and $v_3$; these formulas are used to prune the enumerations of Section~\ref{sec:evidence}, and they show that finite type invariants of degree at most three cannot obstruct $s(K)=3$ (Proposition~\ref{prop:a2v3-realize}).

\subsection{The coefficient \texorpdfstring{$a_2$}{a2} and non-triviality}
Let $\nabla_K(z)=1+a_2(K)z^2+\cdots$ be the Conway polynomial. For a $3$-segment presentation $D$ with the stacking $e_1>e_2>e_3>e_1$, oriented so that it is traversed along $e_1,e_2,e_3$, and for $x\in C_{12}$, $y\in C_{23}$, we write $y\prec_2 x$ if $e_2$ passes through $y$ before $x$. Let $\operatorname{sgn}$ denote the sign of an oriented crossing.

\begin{proposition}\label{prop:a2}
For every $3$-segment presentation $D$ of $K$ which follows the cyclic rule,
\[
a_2(K)=\sum_{\substack{x\in C_{12},\ y\in C_{23}\\ y\prec_2 x}}\operatorname{sgn}(x)\operatorname{sgn}(y).
\]
By cyclic relabelling, $a_2(K)$ is also the analogous sum over $C_{23}\times C_{31}$, ordered along $e_3$, and over $C_{31}\times C_{12}$, ordered along $e_1$.
\end{proposition}

\begin{proof}
For $S\subseteq C_{12}$ and $T\subseteq C_{23}$ let $D_{S,T}$ be obtained from $D$ by changing the crossings of $S\cup T$. By Proposition~\ref{prop:unknotting}, $D_{C_{12},\emptyset}$, $D_{\emptyset,C_{23}}$ and $D_{C_{12},C_{23}}$ are descending (the last one has the stacking $e_3>e_2>e_1$), so they represent the trivial knot and
\[
a_2(K)=a_2(D_{\emptyset,\emptyset})-a_2(D_{C_{12},\emptyset})-a_2(D_{\emptyset,C_{23}})+a_2(D_{C_{12},C_{23}}).
\]
Changing the crossings of $C_{12}$ and of $C_{23}$ one at a time, this mixed difference becomes a sum over $x\in C_{12}$ and $y\in C_{23}$ of mixed differences at the two crossings $x,y$, the other crossings being fixed in some state. By the skein relation $\nabla_{L_+}-\nabla_{L_-}=z\nabla_{L_0}$, applied at $x$ and at $y$, the mixed difference of $\nabla$ at $x,y$ is $\operatorname{sgn}(x)\operatorname{sgn}(y)\,z^2\nabla_{L_{00}}$, where $L_{00}$ is obtained by smoothing $x$ and $y$. Hence its $z^2$-coefficient is $\operatorname{sgn}(x)\operatorname{sgn}(y)$ times the constant term of $\nabla_{L_{00}}$, which is $1$ if $L_{00}$ is a knot and $0$ otherwise. The number of components of $L_{00}$ depends only on the projection: $L_{00}$ is a knot if and only if the two chords of $x$ and $y$ on the Gauss circle cross. The Gauss circle is traversed along $e_1,e_2,e_3$; the chord of $x$ joins a point of $e_1$ to a point of $e_2$, and the chord of $y$ a point of $e_2$ to a point of $e_3$, so they cross exactly when the endpoint of $y$ on $e_2$ lies between the two endpoints of $x$, that is, when $y\prec_2 x$.
\end{proof}

\begin{corollary}\label{cor:a2-non-trivial}
Let $D$ be a $3$-segment presentation of $K$ which follows the cyclic rule.
\begin{enumerate}
\item If the sum of Proposition~\ref{prop:a2} is non-zero, then $K$ is non-trivial.
\item In particular, suppose that all crossings of $C_{12}$ have the same sign and all crossings of $C_{23}$ have the same sign. If $e_2$ passes through some crossing of $C_{23}$ before some crossing of $C_{12}$, then $K$ is non-trivial; indeed $|a_2(K)|$ is the number of such pairs. The same holds for any two of the three sets $C_{12},C_{23},C_{31}$, ordered along the common segment.
\item $|a_2(K)|\le\min\{c_{12}c_{23},\ c_{23}c_{31},\ c_{31}c_{12}\}$.
\end{enumerate}
\end{corollary}

\begin{remark}[Higher coefficients]\label{rem:ckr}
All coefficients of $\nabla_K$ can be read off from the order of the passages, the over/under information and the signs, by the state sum of Chmutov, Khoury and Rossi \cite{CKR}. Choose a base point on the diagram. For a set $S$ of crossings, smooth every crossing of $S$ according to the orientation, and call $S$ \emph{ascending one-component} if the resulting curve is connected and, travelling along it from the base point, each crossing of $S$ is reached first along its over-passage. Then
\[
a_{2m}(K)=\sum_{\substack{|S|=2m\\ S\ \text{ascending one-component}}}\ \prod_{x\in S}\operatorname{sgn}(x)
\]
(for classical knots the same holds with ``under-passage'', \cite{CKR}). For a $3$-segment presentation following the cyclic rule, with base point $p_1$ and $m=1$, the ascending one-component pairs are exactly the pairs $x\in C_{12}$, $y\in C_{23}$ with $y\prec_2x$, which recovers Proposition~\ref{prop:a2}; we checked this equality of sets on $400$ random presentations, and the formula for $a_2,a_4$ on $300$ and for $a_2,a_4,a_6$ on $60$ further random presentations, against the Conway polynomial obtained from knot Floer homology.
For $m\ge2$ the cyclic order of the smoothed curve across all three segments matters. For example, the presentation of $8_{14}$ with $e_1=[0,6]$, $e_2$ visiting $6,4,7,-2,1,3,8$ (first half-circle in the lower half-plane) and $e_3$ visiting $8,5,-1,2,0$ (first half-circle in the upper half-plane), in the half-circle notation of Subsection~\ref{subsec:all-two-bridge}, has $9$ crossings, $(c_{12},c_{23},c_{13})=(3,4,2)$, and no move I or II removes crossings. Numbering the crossings by their first passage from $p_3$, its signed Gauss words are
\[
e_1:\ 1\,\bar2\,3\,4\,\bar5,\qquad e_2:\ 6\,\bar4\,7\,8\,\bar1\,9\,\bar3,\qquad e_3:\ \bar7\,5\,\bar6\,\bar8\,2\,\bar9,
\]
with signs $-,+,+,+,+,+,-,+,-$ at $1,\dots,9$. With base point $p_1$, the eight pairs of Proposition~\ref{prop:a2} contribute $0$ in total, while the ascending one-component sets of four crossings are $\{2,4,6,9\}$, $\{2,3,6,9\}$, $\{2,3,7,9\}$, $\{2,3,8,9\}$, with sign products $-1,-1,+1,-1$. Hence $a_2=0$ and $a_4=-2$, in accordance with $\nabla_{8_{14}}(z)=1-2z^4$, and the presentation is detected as non-trivial by $a_4$ but not by $a_2$. The set $\{2,4,6,9\}$ uses crossings of all three sets $C_{12}$, $C_{23}$, $C_{31}$.
\end{remark}

This is the bound used for pruning in the enumerations of Section~\ref{sec:evidence}. On the other hand, a presentation which cannot be simplified by moves I and II may still represent the trivial knot.

\begin{figure}[ht]
\centering
\includegraphics[width=.62\linewidth]{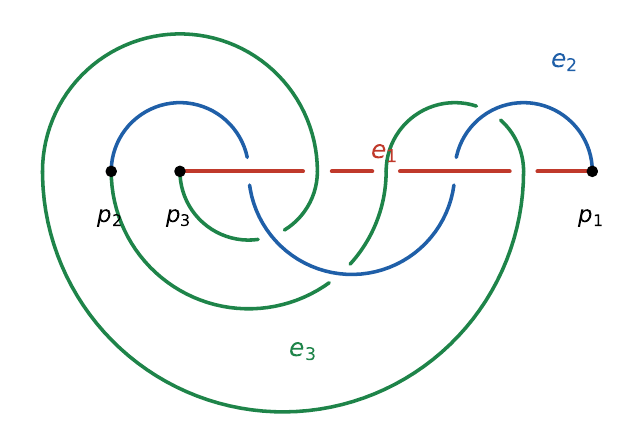}
\caption{An $8$-crossing $3$-segment presentation of the trivial knot with $(c_{12},c_{23},c_{13})=(2,3,3)$. There are no monogons, and the over-strand changes at both bigons. In the half-circle notation of Subsection~\ref{subsec:all-two-bridge}, $e_1=[0,6]$, $e_2$ visits $6,4,1,-1$ with its first half-circle in the upper half-plane, and $e_3$ visits $-1,3,5,-2,2,0$ with its first half-circle in the lower half-plane.}
\label{fig:reduced-unknot}
\end{figure}

\begin{example}\label{ex:reduced-unknot}
The sum in Proposition~\ref{prop:a2} has four terms in Figure~\ref{fig:reduced-unknot}, with signs $-,+,-,+$. The sums along $e_3$ and $e_1$ also vanish. The knot group is infinite cyclic, so the diagram represents the trivial knot. No move I or II decreases its crossing number. The presentation follows the cyclic rule and every segment has crossings, so no edge contains two partition points. Lemma~\ref{lem:rigid} therefore excludes a move III or a partition-point slide that preserves the cyclic rule. Every pair of segments has crossings of both signs. In our experiments, among about $3.5\times10^4$ reduced non-alternating presentations in which the crossings of at least one pair have a common sign, none represented the trivial knot, whereas all of the roughly $170$ reduced presentations of the trivial knot found had crossings of both signs in each pair. A sign-coherent pair consisting of a single crossing does not suffice, however; see Question~\ref{q:coherent} and the example following it.
\end{example}

\subsection{Finite type invariants}
The argument of Proposition~\ref{prop:a2} applies to every finite type invariant. Write $X=C_{12}$, $Y=C_{23}$, $Z=C_{31}$, and for $F\subseteq X\cup Y\cup Z$ and $S\subseteq F$ let $D_{S;F}$ be the singular diagram obtained from $D$ by making the crossings of $S$ double points and changing the crossings of $F\setminus S$. A finite type invariant $v$ is extended to singular knots by $v(D_\times)=v(D_+)-v(D_-)$, and $\varepsilon(x)$ denotes the sign of the crossing $x$ in $D$.

\begin{proposition}\label{prop:finite-type}
Let $D$ be a $3$-segment presentation of $K$ which follows the cyclic rule, and $v$ a finite type invariant of degree at most $d$.
\begin{enumerate}
\item For any two of the sets, say $X$ and $Y$,
\[
v(K)-v(U)=\sum_{\substack{S\subseteq X\cup Y,\ |S|\le d\\ S\cap X\ne\emptyset\ne S\cap Y}}\Bigl(\prod_{x\in S}\varepsilon(x)\Bigr)\,v(D_{S;X\cup Y}),
\]
where $U$ is the trivial knot.
\item If $d\le3$, then
\[
v(K)-v(K^*)=\sum_{x\in X,\ y\in Y,\ z\in Z}\varepsilon(x)\varepsilon(y)\varepsilon(z)\,W_v(G_{x,y,z}),
\]
where $K^*$ is the mirror image, $G_{x,y,z}$ is the chord diagram of the three chords $x,y,z$ on the Gauss circle, and $W_v$ is the weight system of $v$ in degree $3$. Only triples whose chords form a connected intersection graph contribute.
\end{enumerate}
\end{proposition}

\begin{proof}
At each crossing $x$, $v(D)=v(D')+\varepsilon(x)v(D_\times)$, where $D'$ is obtained by changing $x$. Expanding over all crossings of $F$ gives $v(D)=\sum_{S\subseteq F}\bigl(\prod_{x\in S}\varepsilon(x)\bigr)v(D_{S;F})$, and the terms with $|S|>d$ vanish.
Changing all crossings of one or of two of the sets $X,Y,Z$ gives a diagram whose segments are totally ordered by the over/under relation (for instance $e_2>e_3>e_1$ after changing $X$, and $e_3>e_2>e_1$ after changing $X\cup Y$), hence a descending diagram of the trivial knot (Proposition~\ref{prop:unknotting}); changing all crossings gives a diagram of $K^*$.
For (1), the mixed difference $v(D)-v(D^{X})-v(D^{Y})+v(D^{X\cup Y})=v(K)-v(U)$, where $D^T$ is obtained by changing $T$, keeps exactly the terms of the expansion over $F=X\cup Y$ which meet both $X$ and $Y$.
For (2), the mixed difference over the three sets equals $v(K)-3v(U)+3v(U)-v(K^*)$ and keeps the terms meeting $X$, $Y$ and $Z$; for $d\le3$ these are the sets $S=\{x,y,z\}$ with one crossing from each set, and the value of $v$ on a knot with three double points depends only on its chord diagram. A three-chord diagram with disconnected intersection graph has a chord with adjacent endpoints: an isolated chord has this property if the other two interlace, and otherwise the non-interlacing chords include one with adjacent endpoints. The one-term relation gives the asserted vanishing.
\end{proof}

For $d=2$ and $v=a_2$, part~(1) is Proposition~\ref{prop:a2}. For three chords $x\in X$, $y\in Y$, $z\in Z$, the chords of $x$ and $y$ cross if and only if $y\prec_2x$ as before, those of $y$ and $z$ if and only if $z\prec_3y$, and those of $z$ and $x$ if and only if $x\prec_1z$, where $\prec_i$ is the order along $e_i$. Hence a triple is a \emph{triangle} (three pairwise crossing chords) if all three relations hold, and a \emph{path} if exactly two hold. Put
\[
T(D)=\sum_{\substack{x\in X,\,y\in Y,\,z\in Z\\ y\prec_2x,\ z\prec_3y,\ x\prec_1z}}\varepsilon(x)\varepsilon(y)\varepsilon(z),\qquad
P(D)=\sum_{\text{paths}}\varepsilon(x)\varepsilon(y)\varepsilon(z).
\]

\begin{figure}[ht]
\centering
\includegraphics[width=.92\linewidth]{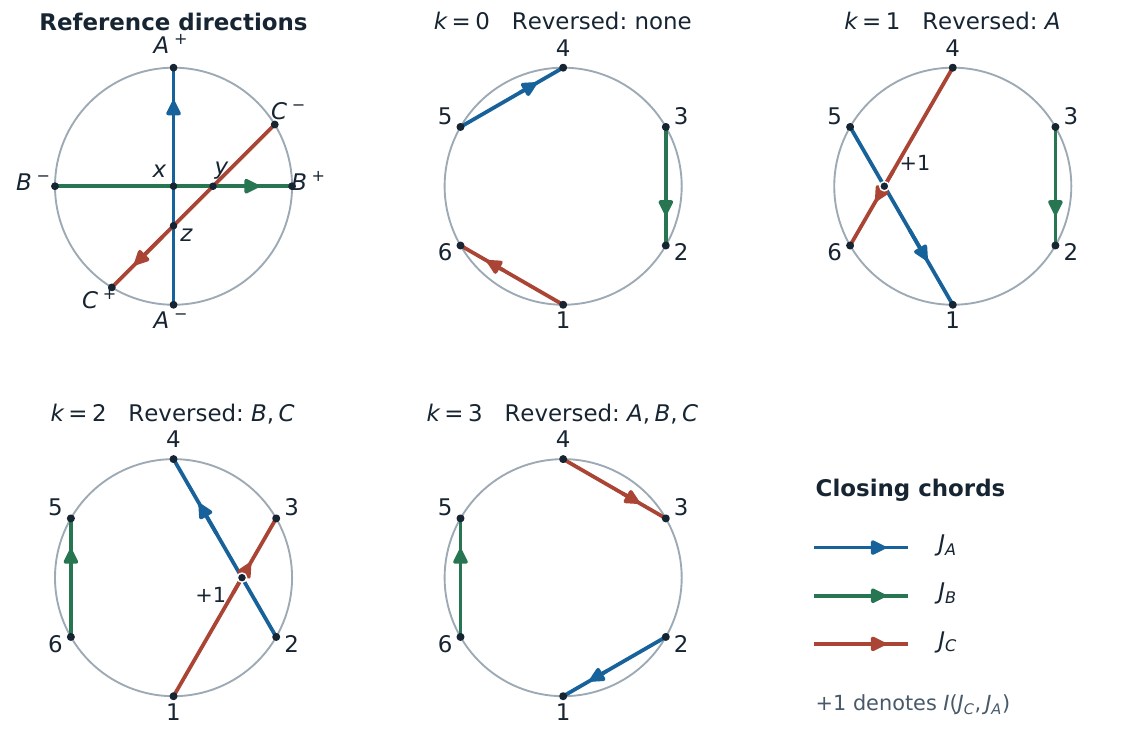}
\caption{The triple point in the proof of Proposition~\ref{prop:path-vanishes}. Top left: the reference directions, in which the local triangle $x\to y\to z\to x$ is clockwise and all three local planar signs are $-1$; the boundary points of the move disk are, counterclockwise from the bottom, $A^-,B^+,C^-,A^+,B^-,C^+$, numbered $1,\dots,6$. The other panels show the closing chords $J_A\colon B^-\to A^+$, $J_B\colon C^-\to B^+$, $J_C\colon A^-\to C^+$ for $k=0,1,2,3$, after exchanging the labels $\pm$ of the reversed strands. Only for $k=1,2$ do two closing chords meet, namely $J_C$ and $J_A$, with $I(J_C,J_A)=+1$; hence $E=0,-1,-1,0$.}
\label{fig:triple-point}
\end{figure}

\begin{proposition}\label{prop:path-vanishes}
For every planar $3$-segment presentation \(D\), the signed path count
\(P(D)\) is zero.
\end{proposition}

\begin{proof}
Use the planar signs \(\operatorname{sgn}\det(\dot e_i,\dot e_{i+1})\)
at intersections of \(e_i,e_{i+1}\), with indices modulo three.
For each pair of segments, the crossing signs differ from the planar signs
by a fixed factor in $\{+1,-1\}$, since the over-strand is fixed for that pair.
Every term in $P(D)$ uses one crossing from each pair, so the product of the
three factors is common to all terms and does not affect vanishing. Regard the segments as
three oriented simple arcs on \(S^2\) forming a closed curve.

We first show that \(P\) is invariant under generic deformations through
such triples of simple arcs. At a tangency of two arc interiors, the
two new crossings have opposite signs and adjacent visits on both
arcs. A path triple uses at most one of them, and substituting one
for the other gives the same order comparisons. A crossing gained at
a common endpoint is last on one arc and first on the next. Its chord
has adjacent ends, so it belongs to no path triple.

Consider a triple-point event involving all three arcs. Write
\(A=e_1,B=e_2,C=e_3\), and \(x=A\cap B,\ y=B\cap C,\ z=C\cap A\)
for the three local crossings. The three order tests
\(a=[y\prec_B x], b=[z\prec_C y], c=[x\prec_A z]\) are all
reversed. All other relative orders are unchanged. Let \(k=a+b+c\)
before the event and let \(s_x,s_y,s_z\) be the planar signs.
The change from the all-local triple, divided by \(s_xs_ys_z\),
is
\[
\delta_k=[k=1]-[k=2].
\]

Take the disk of the move small enough to contain only the three
local strands. Denote the incoming and outgoing endpoints of
each strand on its boundary by \(A^-,A^+\), and similarly for \(B,C\).
In the complementary disk, the oriented parts of the original curve
are three arcs
\[
L_A:A^+\longrightarrow B^-,
\qquad L_B:B^+\longrightarrow C^-,
\qquad L_C:C^+\longrightarrow A^-.
\]
Write \(E=I(L_A,L_B)+I(L_B,L_C)+I(L_C,L_A)\), where \(I\)
is the algebraic intersection number in the given order.

For example, a triple with the two local crossings \(x,y\)
and an exterior crossing \(z'\) changes its path indicator by
\((1-2a)(b'+c'-3b'c')\). Here \(b'\) tests whether \(z'\)
is on the incoming part of \(C\), and \(c'\) whether it is on
the outgoing part of \(A\). The four values of \(b'+c'-3b'c'\)
are the intersection weights in \(E\): \(+1\) for two incoming
parts or two outgoing parts, \(-1\) for incoming \(C\) and
outgoing \(A\), and zero for outgoing \(C\) and incoming \(A\)
(which are parts of the same exterior arc).
The cyclic analogues apply with \(x'\) or \(y'\) exterior.

By reflecting the picture if necessary, assume the local triangle
traversed as \(x\to y\to z\to x\) is clockwise. In the reference
orientations along that triangle every local planar sign is \(-1\).
Comparing the given orientations with the reference gives
\[
s_x=-(1-2c)(1-2a),\qquad
s_y=-(1-2a)(1-2b),\qquad
s_z=-(1-2b)(1-2c).
\]
Thus \((1-2a)/s_z=(1-2b)/s_x=(1-2c)/s_y=(-1)^{k+1}\).
The triples with exactly two local crossings therefore give the
change \((-1)^{k+1}E\), after division by \(s_xs_ys_z\).
Triples with at most one local crossing do not change their contribution.

To evaluate \(E\), close each \(L_i\) by an oriented chord \(J_i\)
in the move disk, from its terminal endpoint to its initial endpoint.
Algebraic intersections of closed curves on the sphere are zero, so
\(E=-I(J_A,J_B)-I(J_B,J_C)-I(J_C,J_A)\).
In the reference orientations the boundary order, counterclockwise, is
\[
 A^-,\ B^+,\ C^-,\ A^+,\ B^-,\ C^+.
\]
Changing one strand orientation exchanges its incoming and outgoing
labels in this list. For \(k=0,3\) the three closing chords are
disjoint. For \(k=1\) take the reversed strand to be \(A\), and
for \(k=2\) take the reversed strands to be \(B,C\), by cyclic
relabelling. In both cases only \(J_C,J_A\) meet, with sign \(+1\)
in that order (Figure~\ref{fig:triple-point}). This follows directly from their alternating boundary
endpoints in the displayed order with the stated exchanges.
Consequently \(E=0,-1,-1,0\) for \(k=0,1,2,3\), respectively.
In every case \(\delta_k+(-1)^{k+1}E=0\). This proves invariance at
the triple point.

Finally, keep the three partition points fixed. Each simple arc
can be isotoped to the corresponding edge of an embedded triangle
with these vertices while avoiding the third partition point:
the complement of that point in \(S^2\) is a plane, in which simple
arcs between two given points are isotopic relative to their ends.
Perform the three isotopies successively. A generic deformation has
only the tangencies, triple points and events at common endpoints
just considered; each segment stays simple and avoids the third
partition point. The terminal diagram has no crossings, hence
\(P=0\).
\end{proof}

\begin{corollary}\label{cor:v3}
Let $v_3(K)$ be the coefficient of $h^3$ in $V_K(e^h)$, where $V_K$ is the Jones polynomial. For every $3$-segment presentation $D$ of $K$ which follows the cyclic rule,
\[
v_3(K)=-6\,T(D).
\]
In particular, if $T(D)\ne0$, then $K$ is not amphichiral, and hence non-trivial.
\end{corollary}

\begin{proof}
The invariant $v_3$ has degree $3$ and $v_3(K^*)=-v_3(K)$, since $V_{K^*}(t)=V_K(t^{-1})$. Its weight system takes the values $-12$ on the triangle and $-6$ on the path; we computed these values directly as alternating sums of $v_3$ over the eight resolutions of three double points, and they satisfy the four-term relation, by which the path is half the triangle in degree $3$. Proposition~\ref{prop:finite-type}(2) gives $2v_3(K)=-12\,T(D)-6\,P(D)$, and Proposition~\ref{prop:path-vanishes} gives the claim.
\end{proof}

For the minimal diagram of the right-handed trefoil, which is a $3$-segment presentation with one crossing in each set, all positive and pairwise crossing chords, $T=1$ and $P=0$, and indeed $v_3=-6$ for $V(t)=t+t^3-t^4$. We checked Corollary~\ref{cor:v3} against the Jones polynomial on $450$ random presentations.

For every $3$-segment presentation following the cyclic rule, Corollary~\ref{cor:v3} gives $|v_3(K)|/6=|T(D)|\le c_{12}c_{23}c_{31}\le (c(D)/3)^3$, and Proposition~\ref{prop:a2} gives $|a_2(K)|\le\min\{c_{12}c_{23},c_{23}c_{31},c_{31}c_{12}\}$. Such bounds prune searches with a bounded number of crossings, but the values of $a_2$ and $v_3$ alone cannot show $s(K)>3$.

\begin{proposition}\label{prop:a2v3-realize}
Put $w(K)=-v_3(K)/6$, so that $w$ is integer valued and $w=1$ for the right-handed trefoil. For every pair of integers $(m,n)$ there is a non-trivial two-bridge knot $K$, hence a knot with $s(K)=3$, such that $a_2(K)=m$ and $w(K)=n$.
\end{proposition}

\begin{proof}
For integers $x,m$ let $K_{x,m}=C(2x,2,-2x,2x,2,-2x-2m)$. The formula of Ichihara and Wu for the order-$3$ invariant of two-bridge knots of genus $3$ \cite[Section~4]{IW}, together with the corresponding formula for $a_2$, gives
\[
a_2(K_{x,m})=m,\qquad w(K_{x,m})=-x-\tfrac{m(m+1)}2
\]
up to a global sign of $w$ which depends only on the convention for $C(\cdots)$ and is removed by taking mirror images; for $m=0$ these are the knots $C(2x,2,-2x,2x,2,-2x)$ of \cite{IS}. Given $(m,n)$, choose $x=-n-m(m+1)/2$. If $x\ne0$ and $x+m\ne0$, the knot $K_{x,m}$, or its mirror image, has the required values. If $x+m=0$, the last coefficient vanishes and $K_{x,m}=C(2x,2,-2x,2x)$; if $x=0$, the expansion reduces to the twist knot $C(-2m,2)$ or its mirror image; in both cases the same values result, and $(0,0)$ is realised by the non-trivial knot $C(2,2,-4,2)=S(31,12)$. By Theorem~\ref{thm:twobridge-intro} all these knots have segment number $3$.
As an independent check, we computed $a_2$ by Proposition~\ref{prop:a2} and $w=T$ by Corollary~\ref{cor:v3} directly on the signed double spirals of Theorem~\ref{thm:all-two-bridge} for $|x|,|m|\le3$ and for the reduced cases; with our conventions $w$ has the opposite sign, uniformly.
\end{proof}

Thus an obstruction to $s(K)=3$ cannot come from the finite type invariants of degree at most three alone; it has to use other information, or a bound which does not degenerate as the number of crossings grows.

\section{Experiments}\label{sec:evidence}

\subsection{Random searches}
By Propositions~\ref{prop:curves} and \ref{prop:normal}, the knots with $s(K)\le 3$ are, up to mirror image, the knots $K_+(C)$ of tri-simple curves $C$.
We generated $82{,}985$ tri-simple curves by choosing the three segments as random self-avoiding polygonal arcs between three fixed points, with $6$ to $14$ edges each (up to $105$ crossings); $|\sigma|$ took the values $0,2,4$ with frequencies $56{,}225$, $26{,}679$, $81$, and $|\tau|$, computed for $10{,}558$ of them, the values $0,1,2$.
Random walks in the space of tri-simple curves (moving, inserting and deleting vertices, rejecting moves that create a self-crossing of a segment) which maximise $|\sigma|$ or $\col_3$ instead readily reach $|\sigma|=6$ (Example~\ref{ex:sig}) and $\col_3=3^4$ (Proposition~\ref{prop:K14}).

\subsection{Census of small knots}
A further random walk recorded the knot type of every diagram with at most ten crossings after simplification, together with the smallest presentation found, and a simulated annealing search was directed towards specific knots using the hyperbolic volume and the determinant.
Finally, starting from a minimal diagram of a given knot, we applied random Reidemeister moves (with the \texttt{backtrack} and \texttt{simplify} routines of spherogram) and minimised the number of crossings violating the conditions of Definition~\ref{def:segment}, over all choices of the three partition points and of the cyclic rule; this produced the presentations of $8_3$ ($14$ crossings) and of $5_1$ with $9$ crossings.
The presentation of $7_3$ ($14$ crossings) and the first presentations of $8_5$ ($33$ crossings) and $8_{10}$ ($27$) were obtained by re-routing passes, starting from the standard diagrams: an arc which passes over (or under) everything it meets is replaced by another such arc with the same end points, chosen so as to decrease the number of crossings violating Definition~\ref{def:segment}; the optimal re-routing of a single pass is a shortest-path problem in the dual graph of the remaining diagram. Read backwards, each of these sequences is a certificate of the knot type (Section~\ref{sec:method}). The presentations of $8_5$ and $8_{10}$ were later replaced by presentations with $23$ and $21$ crossings found by the targeted enumeration below.
The results are in Table~\ref{tab:census}; per crossing number, 1 of 1 with $c=3$, 1 of 1 with $c=4$, 2 of 2 with $c=5$, 3 of 3 with $c=6$, 7 of 7 with $c=7$, 21 of 21 with $c=8$, 48 of 49 with $c=9$, 127 of 165 with $c=10$.
The presentations of $9_2$ and $10_{10}$ were obtained as corner twists (Proposition~\ref{prop:corner}) of those of $3_1$ and $6_3$.
In particular every prime knot with at most eight crossings has segment number three (Figures~\ref{fig:gallery} and~\ref{fig:gallery-large}). The last two cases, $8_{15}$ and $8_{19}=T(3,4)$, two Montesinos knots of length three with $u(8_{15})=2$ and $u(8_{19})=3$, needed much larger presentations, with $27$ and $30$ crossings.

The found knots include knots with long twist regions ($8_1=C(6,2)$, $8_2=C(5,1,2)$, $8_3=C(4,4)$), the torus knots $5_1$, $7_1$ and $9_1$, and the length-three Montesinos knots $8_5$, $8_{10}$, $8_{20}$ and $8_{21}$, so none of these features alone accounts for the difficulty of the last two cases.
For $8_{15}$ and $8_{19}$ the first re-routing searches, with up to $60$ crossings, reached only diagrams with a single violating crossing; in almost all of these the violating crossing joins two different segments, and changing it gives a $3$-segment presentation of a simpler knot ($5_1$ or $5_2$). The presentations of $8_{15}$ ($27$ crossings, $(c_{12},c_{23},c_{13})=(10,5,12)$) and of $8_{19}$ ($30$ crossings, $(5,7,18)$) were finally found by a combination of these re-routing searches with the exhaustive enumeration below; they are certified by $8$ and $11$ re-routings of passes, respectively, ending at the standard diagrams of the Knot Atlas without mirror image (Section~\ref{sec:method}).

\subsection{Exhaustive enumeration}
Let $D$ be a $3$-segment presentation with a minimal number of crossings for its knot type. Then $D$ has no monogon, and no bigon whose two edges belong to different segments and contain no partition point, since a Reidemeister move~I or~II would give a smaller $3$-segment presentation (by the cyclic rule, the same segment is the over-strand at both crossings of such a bigon); we call such diagrams \emph{reduced}. If $u(K)\ge2$, then every pair of segments of $D$ crosses at least twice (Corollary~\ref{cor:c3}).
We enumerated all reduced $3$-segment presentations with at most $17$ crossings in which every pair of segments crosses at least twice, as follows. The arc $e_1$ is fixed in the sphere, and $e_2$ is grown from the end of $e_1$: each time it is to cross $e_1$, we choose, by walking along the boundary of the face containing its tip, which edge of $e_1$ it crosses; since an arc in a disc is determined up to isotopy by its end points, every diagram is obtained. After $c_{12}$ crossings the partition point is placed and $e_3$ is grown in the same way, crossing $e_1$ and $e_2$ the prescribed numbers of times, and the diagram is closed up if the tip lies in the face containing the start of $e_1$. Cyclic relabelling, reversal of orientation and mirror image allow us to assume $c_{12}\le c_{23}$ and $c_{12}\le c_{13}$.
The enumeration produced $106\,322$ reduced diagrams with $6$ to $17$ crossings, every one of which was identified with SnapPy; a diagram replaces the entry of Table~\ref{tab:census} when it has fewer crossings than the presentation found before (the earlier presentations remain in the supplementary data).
As a check of completeness, every presentation of the earlier census with at most $14$ crossings and with all $c_{ij}\ge2$ was found again, with at most as many crossings. As an independent check, we generated all tri-simple curves with at most $10$ crossings (there are $1,4,13,47,165,637,2433,9602,37959,151904$ of them with $1,\dots,10$ crossings, up to homeomorphism of the sphere and relabelling) by the inverse of monotone simplification: starting from a simple closed curve with three marked points, we inserted empty monogons at partition points and empty bigons, and closed up under triangle moves and slides of partition points across crossings, keeping the three arcs simple. The reduced curves with all $c_{ij}\ge2$ among them coincide, as sets of equivalence classes, with those of the enumeration ($2,4,4,32,55$ classes with $6,\dots,10$ crossings). Moreover, each of the $523$ classes of reduced curves with all $c_{ij}\ge 2$ and at most $12$ crossings can be brought, by at most one triangle move or slide, to a curve with an empty monogon or bigon.
For a knot $K$ with $u(K)\ge2$, a presentation with $c_3(K)$ crossings is reduced and has all $c_{ij}\ge2$. If $c_3(K)\le17$, the enumeration therefore contains a minimal presentation, whose identification is subject to the checks discussed below.
To guard against misidentification, for each such knot $K$ every enumerated diagram with fewer crossings whose determinant and knot Floer homology (total rank, genus, fibredness) agree with those of $K$ was compared with the exterior of $K$ using SnapPy; there were $1096$ such diagrams. The limits of negative isometry tests are discussed below.

\begin{theorem}[computer-assisted]\label{thm:enum}
\begin{enumerate}
\item The $105$ prime knots of Table~\ref{tab:census} found by the enumerations described in this section have segment number $3$: $55$ knots found by the enumeration up to $17$ crossings, $35$ knots found with $18$ to $22$ crossings and $15$ knots found with $23$ to $26$ crossings by the extended searches described below. They are listed in the supplementary data (\texttt{enumeration/theorem\_8\_1\_lists.txt}).
\item If the entry $N$ in Table~\ref{tab:census} equals the crossing number $c(K)$, then $c_3(K)=N$. In particular, $c_3(8_{18})=8$.
\item $c_3(8_{15})=27$ and $29\le c_3(8_{19})\le30$.
\end{enumerate}
\end{theorem}

For~(3), the upper bounds are given by the presentations in Table~\ref{tab:census}. For the lower bounds, $u(8_{15})=2$ and $u(8_{19})=3$, so all pairs of segments of a minimal presentation cross at least $2$, respectively $3$, times. For $8_{19}$, no diagram with at most $21$ crossings has determinant $3$ and the knot Floer homology of $T(3,4)$ (total rank $5$, genus $3$, fibred); for $8_{15}$, no diagram with at most $17$ crossings of the enumeration above is $8_{15}$. The remaining crossing counts, $18$ to $26$ for $8_{15}$ and $22$ to $28$ for $8_{19}$, are covered by the filtered searches with Jones polynomial exclusions described next. All $4{,}137$ archived candidates of the filtered searches for these two knots, including those of earlier ranges ($4{,}071$ of them lie in the ranges $18$--$26$ and $22$--$28$ above), were re-checked with the Jones polynomial, and none of them is $8_{15}$ or $8_{19}$. The enumeration was continued, with filters, to $26$ crossings for $8_{15}$ and to $28$ crossings for $8_{19}$: only diagrams with the determinant and the value $\Delta(2)\Delta(1/2)$ of the target knot were kept, branches were pruned using the fact that $a_2=v_2$ is a signed count of pairs of crossings of $e_1$ with $e_2$ and of $e_2$ with $e_3$ in a prescribed relative position, which can be bounded during the construction, and all remaining candidates were excluded by values of the Jones polynomial; up to $23$ crossings the $550$ candidates for $8_{15}$ were also excluded independently by SnapPy.
The enumeration was also continued to $22$ crossings, keeping only diagrams whose fingerprint $(|\Delta(-1)|,\Delta(2)\Delta(1/2))$, computed from the Alexander matrix, agrees with that of a prime knot with at most ten crossings not yet in the census; the $75{,}214$ diagrams with a matching fingerprint among the $7{,}685{,}052$ reduced diagrams with $18$ to $22$ crossings were identified with SnapPy. This gave the $35$ knots of~(1) with $18$ to $22$ crossings; each of these presentations was checked independently as described in Section~\ref{sec:method}.
Finally, the fingerprint test was moved into the enumeration program itself, with the fingerprint reduced modulo $1009$ (up to sign, since the Alexander polynomial is defined only up to $\pm t^k$), for the $54$ prime knots with at most ten crossings still missing; this made it possible to enumerate all reduced diagrams with all $c_{ij}\ge2$ and $23$ to $26$ crossings (for $26$ crossings and $\min c_{ij}=8$ with the pruning described below). The $240{,}776$ diagrams that passed the test were identified with SnapPy, giving the $15$ knots of~(1) with $23$ to $26$ crossings; they are checked in the same way as the others.
The pruning by $a_2$ used for $8_{19}$ makes such searches much faster once a target knot is fixed, since the signed count $a_2$ is bounded during the construction of $e_3$. With this pruning we also enumerated all reduced diagrams with $\min c_{ij}\ge u$ and at most $26$ crossings with the determinant, the value $\Delta(2)\Delta(1/2)$ and the coefficient $a_2$ of $9_1$, and at most $30$ crossings with those of $T(2,11)$; all candidates were identified with SnapPy, and the first presentations of $9_1$ appear at $27$ crossings, none of $T(2,11)$ at $30$ or fewer. These searches support $c_3(9_1)=27$ and $c_3(T(2,11))\ge31$, subject to the identification caveat below.
\medskip\noindent\textbf{Targeted enumeration.}
To test minimality for the remaining entries of Table~\ref{tab:census}, we first removed all monogons and all bigons whose two crossings have the same over-strand from the presentations of the table; these are the crossing-decreasing Reidemeister moves. This reduced $4_1$ from $5$ to $4$ crossings, $5_2$ from $7$ to $6$ and $10_{92}$ from $39$ to $19$.
Then, for every knot $K$ of the table with $c_3(K)$ not yet determined and with $N(K)$ crossings in the table, we ran the enumeration with the pruning by $a_2$ and the filter by the fingerprint $(|\Delta_K(-1)|,\Delta_K(2)\Delta_K(1/2))$ modulo $1009$ for all numbers of crossings $N$ with $N_0(K)\le N<\min(N(K),26)$, where $N_0(K)=18$ if $u(K)\ge2$ (by the enumeration above) and $N_0(K)=c(K)$ if $u(K)=1$; in the latter case the program was extended to allow $\min c_{ij}=1$. The searches found presentations of $8_5$ with $23$ crossings and of $8_{10}$ with $21$ crossings, and no presentation of any other target with fewer crossings than in the table; for $10_{57}$, $10_{76}$, $10_{109}$ and $10_{125}$, which have $26$ crossings in the table, the search covered $18$ to $25$ crossings. The $9{,}469$ candidates were identified with SnapPy; independently, their knot Floer homology (total rank, genus, fibredness, $|\tau|$) was compared with that of the targets, and candidates with matching knot Floer data were compared further. Some simplify to a diagram with fewer crossings than the crossing number of the target (for instance the square knot and $8_{20}$ when the target is $10_{140}$). Other comparisons used SnapPy's isometry test.
As positive controls, the searches recover presentations of $10_{92}$ with $19$, of $8_7$ with $12$ and of $9_{43}$ with $23$ crossings. As an independent check of the extension to $\min c_{ij}=1$, a direct enumeration of signed Gauss words with $\min c_{ij}=1$ and at most $9$ crossings, without pruning and with planarity tested by the algorithm of de~Fraysseix and Ossona de Mendez \cite{dFOM}, gives the same first appearances of the knots concerned ($4_1$ at $4$, $5_2$ at $6$, $6_1$ at $7$ and $7_2$ at $9$ crossings).
Taken together, the searches support the bold entries of Table~\ref{tab:census} as minima. The exact claims in Theorem~\ref{thm:enum} are limited to the cases stated there.

\medskip\noindent\textbf{Completeness of the pruning and status of the exclusions.}
The searches use the following pruning tests: (i) the diagram contains an empty monogon, or a bigon whose two crossings have the same over-strand and whose sides contain no partition point; (ii) the prescribed count for a pair falls below the search threshold, in particular $u(K)$ when the enumeration uses the unknotting number of the target $K$; (iii) the signed count of Proposition~\ref{prop:a2}, bounded during the construction of $e_3$, can no longer reach $a_2(K)$; (iv) the fingerprint of the closed diagram differs from that of $K$. A presentation of $K$ with $c_3(K)$ crossings follows the cyclic rule, is reduced and has all $c_{ij}\ge u(K)$ (Corollary~\ref{cor:c3}), and it satisfies (iii) and (iv) for $K$; so none of these rules discards it.
A surviving candidate can be ruled out by a distinguishing invariant (for example, the determinant, $\Delta(2)\Delta(1/2)$, knot Floer homology or the Jones polynomial), a positive identification as a different knot up to mirror image, or a simplification to fewer than $c(K)$ crossings. Some comparisons used negative reports from SnapPy's numerical isometry test, which can miss an existing isometry. Lower bounds relying on those reports remain numerical evidence. The supplementary data contain the presentations, verification scripts and available search logs. A complete record of the exclusion method for each candidate is unavailable for the census as a whole.
Theorem~\ref{thm:enum} supersedes the smaller exhaustive searches reported in the first version of this paper.

\medskip\noindent\textbf{The open cases.}
The $39$ prime knots with at most ten crossings which are not in Table~\ref{tab:census} are
$9_{16}$, $10_{46}$--$10_{51}$, $10_{53}$--$10_{56}$, $10_{61}$--$10_{64}$, $10_{66}$, $10_{78}$--$10_{81}$, $10_{96}$, $10_{98}$--$10_{101}$, $10_{111}$, $10_{120}$, $10_{124}$, $10_{126}$--$10_{128}$, $10_{134}$, $10_{139}$, $10_{142}$, $10_{152}$--$10_{154}$, $10_{162}$, $10_{165}$.
All of them have bridge number three. The enumerations found no reduced $3$-segment presentation of these knots with all $c_{ij}\ge2$ and at most $25$ crossings; for $10_{165}$, which was not among the targets of the filtered searches (see Section~\ref{sec:method} on the names of knots), this applies only to at most $17$ crossings. For the targets with $u(K)\ge2$ other than $10_{165}$, this is evidence for $c_3(K)\ge26$, subject to the identification caveat above.

\section{Certification}\label{sec:method}

Polygonal presentations in this paper are specified by integer vertex coordinates; the remaining presentations are given by PD codes.
For each polygonal presentation a script, available with the supplementary data, checks in exact rational arithmetic that the polygon is generic (no vertex lies on a non-incident edge, no two edges overlap, all double points are transverse and distinct) and that no segment has a self-crossing.
It then lists the crossings, determines the over-strand by the cyclic rule, and writes a PD code, from which Fox colourings are computed by Gaussian elimination over $\mathbb{Z}/p$.
Knot types were identified with SnapPy \cite{SnapPy} (hyperbolic structure and the Hoste--Thistlethwaite census for hyperbolic knots; fundamental group and knot Floer homology for torus knots).
Presentations from the exhaustive enumeration are given by signed Gauss words; a script checks that each word is realised by a diagram in the sphere (by the Euler characteristic), that the three segments are simple and follow the cyclic rule, and identifies the knot with SnapPy. For the coordinate presentations of $30$ of the $35$ prime knots with at most eight crossings, which remain in the supplementary data, the identification does not depend on SnapPy: for each of them we give an explicit sequence of Reidemeister moves, each performed across a face (a monogon, a bigon whose two crossings have the same over-strand, or a non-cyclic triangle), from the $3$-segment presentation to the diagram of the knot in KnotInfo \cite{KnotInfo} or its mirror image.
The certificates contain between $0$ and $36$ transitions (Reidemeister moves or harmless reversals of the Gauss-word orientation); for $3_1$ and $7_7$ the presentation is itself the standard diagram.
Each sequence, given by signed Gauss codes, is checked by the same short program as the certificates of Section~\ref{sec:square}, which uses no knot-theoretic software.
For the remaining five, $7_3$, $8_5$, $8_{10}$, $8_{15}$ and $8_{19}$ (for $8_5$ and $8_{10}$ their earlier presentations with $33$ and $27$ crossings, which remain in the supplementary data, while the present ones, with $23$ and $21$ crossings, come from the enumeration and are checked as the other Gauss-word presentations; for the last two with $8$ and $11$ re-routings, ending at the diagrams of the Knot Atlas without mirror image), the certificate is a sequence of re-routings of passes, as for the granny knot in Section~\ref{sec:square}, from the $3$-segment presentation to the diagram of the knot in the Knot Atlas: at each step an arc passing over (or under) everything it meets is replaced by another such arc with the same end points, and the program checks that the signed Gauss words of the complementary interval agree before and after the step. SnapPy identifies these knots independently ($8_{19}$ by its fundamental group $\langle a,b\mid a^3=b^4\rangle$ and its knot Floer homology).
Knots are named as in KnotInfo \cite{KnotInfo}. SnapPy names knots with at most ten crossings after Rolfsen's table, which differs from KnotInfo in two places: the names $10_{83}$ and $10_{86}$ are interchanged, and SnapPy's $10_{163},\dots,10_{166}$ are KnotInfo's $10_{162},\dots,10_{165}$ (Rolfsen's $10_{161}$ and $10_{162}$ are the same knot). All identifications were converted to KnotInfo names; in this conversion the Alexander polynomial of every presentation of Table~\ref{tab:census}, computed from its knot Floer homology, agrees with that of the knot of KnotInfo of the same name. Since the targeted searches used the invariants of the knot of KnotInfo with the name recorded in the logs, the entries $10_{163}$ and $10_{164}$ (Rolfsen's $10_{164}$ and $10_{165}$) were checked again with the correct targets. For $10_{163}$, with $u=2$, the value $c_3=11$ follows from the complete enumeration up to $17$ crossings, whose identifications are consistent among themselves; directly, no reduced $3$-segment presentation with $10$ crossings has the Alexander polynomial of $10_{163}$. For $10_{164}$, with $u=1$, we enumerated all reduced $3$-segment presentations with $10$ to $13$ crossings, now allowing $c_{ij}=1$ ($3{,}794$ diagrams): the $34$ diagrams with the Alexander polynomial of $10_{164}$ are identified by SnapPy, by an isometry it finds, as $10_{10}$ ($28$ diagrams) and \texttt{K12n333} ($6$ diagrams), so none of them is $10_{164}$; hence $c_3(10_{164})=14$.
Signatures were computed by the Gordon--Litherland formula \cite{GLi} from a checkerboard surface; the implementation was tested against the standard values for all prime knots with at most ten crossings, and on all diagrams of the census.

\section{Questions}\label{sec:questions}

\begin{question}\label{q:main}
Prove or disprove Conjecture~\ref{conj:main}. In particular:
\begin{enumerate}
\item Do the $39$ prime knots with at most ten crossings listed in Section~\ref{sec:evidence}, all of bridge number three, have segment number $3$? Is $s(K)=3$ for every knot of bridge number three?
\item Is $s(T(p,q))=3$ for all torus knots? The first open case is $T(3,5)=10_{124}$.
\item Is $s(K_1\# K_2)=3$ whenever $s(K_1)=s(K_2)=3$? Example~\ref{ex:square} answers this for two trefoils.
\end{enumerate}
\end{question}

\begin{question}
Is there a polynomial $P$ such that $c_3(K)\le P(c(K))$ for every knot $K$ with $s(K)=3$? On twist knots $c_3$ grows at most linearly, and on the knots $T(2,q)$ at most quadratically (Section~\ref{sec:T2q}).
\end{question}

\begin{question}\label{q:T2q}
Is $c_3(T(2,2k+1))=k^2+3k-1$ for all $k\ge1$? This holds for $k=1$, and the searches support the equality for $2\le k\le4$. Is $c_3(T(2,q))/q$ unbounded? Is $c_3(8_{19})=29$ or $30$?
\end{question}

\begin{question}\label{q:coherent}
Let $D$ be a $3$-segment presentation following the cyclic rule in which no Reidemeister move I or II removes crossings, and suppose that one of $C_{12},C_{23},C_{31}$ consists of at least two crossings, all of the same sign. Is $K(D)$ non-trivial? (By Corollary~\ref{cor:a2-non-trivial}, this holds if two of the three sets are sign-coherent and their crossings interleave along the common segment.)
\end{question}

The assumption of at least two crossings cannot be dropped. The presentation with signed Gauss words
\[
e_1:\ \bar1\,\bar2\,3\,\bar4\,\bar5\,\bar6\,\bar7,\qquad
e_2:\ 8\,9\,10\,11\,\bar3\,12,\qquad
e_3:\ 5\,\overline{10}\,1\,\bar8\,7\,2\,\overline{11}\,4\,\overline{12}\,6\,\bar9,
\]
where $k$ and $\bar k$ denote the over- and under-passage through the crossing $k$, and the rotation at the crossings $2,3,5,6,9,10,12$ is (incoming under, incoming over, outgoing under, outgoing over) counterclockwise and the opposite at the others, is a planar $3$-segment presentation with $12$ crossings and $C_{12}=\{3\}$. It has no monogon, and its only bigon, bounded by the crossings $7$ and $8$, cannot be removed by a move II; nevertheless it represents the trivial knot (moves III across the triangles $\{1,9,10\}$ and $\{1,6,7\}$ followed by moves II and I reduce it to a circle). The data are in the supplementary material.

\subsection*{Data availability}
The diagram data (integer coordinates, PD codes or signed Gauss words), verification scripts and available search logs are supplied as ancillary files (directory \texttt{anc}).

\subsection*{Acknowledgements}
The author thanks Takuji Nakamura and the other participants of the Topology Seminar at Tokyo Woman's Christian University, where this work was presented in November 2026; their question whether two-bridge knots, as well as torus knots, have segment number three led to Theorem~\ref{thm:twobridge-intro}. The author also recalls the advice, handed down in the knot theory community and attributed to Kunio Murasugi, that two-bridge knots and torus knots are the two fundamental classes of knots and that a problem should first be attacked for these two classes; it shaped the direction of this work.

The computations use SnapPy, spherogram, the knot Floer homology program of Z.~Szab\'o, and the KnotInfo database \cite{KnotInfo}.
The computer searches, the verification scripts and a first draft of this paper were prepared with the assistance of an AI system (Claude, Anthropic).
The $3$-segment presentations of the granny knot (Figure~\ref{fig:granny}), of the torus knots $7_1$ and $9_1$ and of the knots $7_3$, $8_{15}$ and $8_{19}$ and the first presentations of $8_5$ and $8_{10}$, together with their certificates by re-routing of passes or by Reidemeister moves, the observation of Proposition~\ref{prop:braidlike}, as well as the enumeration code of Section~\ref{sec:evidence} and its extension to $26$ and $28$ crossings for $8_{15}$ and $8_{19}$, were obtained with the assistance of ChatGPT (OpenAI).
The presentations $D_k$ of $T(2,2k+1)$ in Proposition~\ref{prop:T2q}, as well as Lemma~\ref{lem:rigid}, the present presentations of $8_5$, $8_{10}$ and $10_{92}$ and the targeted enumeration of Section~\ref{sec:evidence}, were found with the assistance of Claude (Anthropic).
The signed double spiral of Theorem~\ref{thm:all-two-bridge} and its computer verification were found with the assistance of Claude (Anthropic), and the written proof of Theorem~\ref{thm:all-two-bridge} was prepared with the assistance of ChatGPT (OpenAI).
The idea of building infinite families of $3$-segment plat presentations of every bridge number, leading to Theorem~\ref{thm:every-bridge}, is due to the author; the $6$-plat $\beta_0^{(k)}$, the braid $G$ and the computations in the proof of Theorem~\ref{thm:every-bridge} were obtained with the assistance of ChatGPT (OpenAI), and they were verified independently, and Figure~\ref{fig:K40} drawn, with the assistance of Claude (Anthropic).
The drawings of Figures~\ref{fig:gallery} and~\ref{fig:gallery-large} were produced by circle packing and a planarity-preserving relaxation with the assistance of Claude (Anthropic).
All mathematical statements, proofs and computational results were checked by the author, and the author takes full responsibility for the content.

\appendix
\section{Coordinates, census and gallery}\label{app:coords}

\begin{table}[ht]\centering\small
\caption{Vertices of the polygonal segments of Figure~\ref{fig:K14}, listed in the direction of travel ($e_1$ from $p_3$ to $p_1$, $e_2$ from $p_1$ to $p_2$, $e_3$ from $p_2$ to $p_3$).}\label{tab:coordsK14}
\begin{tabular}{ll}\toprule
$e_1$ & \parbox[t]{0.82\textwidth}{\raggedright (989,-728), (802,-256), (489,177), (434,225), (250,293), (-783,-139), (1357,-964), (889,190), (85,543), (-452,783)}\\[2pt]
$e_2$ & \parbox[t]{0.82\textwidth}{\raggedright (-452,783), (66,718), (712,693), (919,-888), (-407,-373), (436,158), (-99,403), (-725,857), (-871,-414), (-625,-1149)}\\[2pt]
$e_3$ & \parbox[t]{0.82\textwidth}{\raggedright (-625,-1149), (37,-233), (201,377), (-169,537), (-1,433), (-201,495), (-690,837), (569,554), (437,-154), (989,-728)}\\[2pt]
\bottomrule\end{tabular}\end{table}

\begin{table}[ht]\centering\small
\caption{Vertices of the polygonal segments of the $3$-segment presentation of the square knot (Figure~\ref{fig:square}), with the same conventions as Table~\ref{tab:coordsK14}.}\label{tab:coordsSq}
\begin{tabular}{ll}\toprule
$e_1$ & \parbox[t]{0.82\textwidth}{\raggedright (3732,-2076), (1210,-440), (392,1295), (-779,-649), (-775,-156), (-1466,1500)}\\[2pt]
$e_2$ & \parbox[t]{0.82\textwidth}{\raggedright (-1466,1500), (1289,-310), (2283,-939), (771,-702), (-774,-344), (-599,-527), (209,-1843)}\\[2pt]
$e_3$ & \parbox[t]{0.82\textwidth}{\raggedright (209,-1843), (632,-980), (1268,74), (-231,-25), (-819,-786), (-787,-622), (-148,239), (-2166,1465), (-202,1785), (3558,424), (3732,-2076)}\\[2pt]
\bottomrule\end{tabular}\end{table}
\clearpage
{\small\setlength{\tabcolsep}{3pt}
\begin{longtable}{lcccc|lcccc|lcccc}
\caption{Prime knots for which a $3$-segment presentation was found. Here $N$ is the number of crossings of the presentation, and $c_{ij}$ the number of crossings between $e_i$ and $e_j$. Presentations marked $^\dagger$ come from the exhaustive enumeration of Section~\ref{sec:evidence} or from reduction by Reidemeister moves I and II, and are given by signed Gauss words, those marked $^*$ by PD codes, all others by integer coordinates; all are in the supplementary data. A bold $N$ marks an entry supported as minimal by the searches. Theorem~\ref{thm:enum} gives the exact claims, and Section~\ref{sec:evidence} describes the limits of the numerical exclusions.}\label{tab:census}\\
\hline knot & $N$ & $c_{12}$ & $c_{23}$ & $c_{13}$ & knot & $N$ & $c_{12}$ & $c_{23}$ & $c_{13}$ & knot & $N$ & $c_{12}$ & $c_{23}$ & $c_{13}$\\ \hline\endfirsthead
\hline knot & $N$ & $c_{12}$ & $c_{23}$ & $c_{13}$ & knot & $N$ & $c_{12}$ & $c_{23}$ & $c_{13}$ & knot & $N$ & $c_{12}$ & $c_{23}$ & $c_{13}$\\ \hline\endhead
\hline\endfoot
\hline\endlastfoot
$3_{1}$ & $\mathbf{3}$ & 1 & 1 & 1 & $9_{3}^\dagger$ & $\mathbf{24}$ & 3 & 4 & 17 & $9_{41}^\dagger$ & $\mathbf{9}$ & 3 & 3 & 3\\
$4_{1}^\dagger$ & $\mathbf{4}$ & 1 & 1 & 2 & $9_{4}^\dagger$ & $\mathbf{19}$ & 2 & 4 & 13 & $9_{42}^\dagger$ & $\mathbf{12}$ & 3 & 4 & 5\\
$5_{1}^*$ & $\mathbf{9}$ & 3 & 2 & 4 & $9_{5}^\dagger$ & $\mathbf{16}$ & 2 & 3 & 11 & $9_{43}^\dagger$ & $\mathbf{23}$ & 5 & 8 & 10\\
$5_{2}^\dagger$ & $\mathbf{6}$ & 3 & 2 & 1 & $9_{6}^\dagger$ & $\mathbf{20}$ & 3 & 4 & 13 & $9_{44}^\dagger$ & $\mathbf{11}$ & 3 & 3 & 5\\
$6_{1}^\dagger$ & $\mathbf{7}$ & 2 & 2 & 3 & $9_{7}^\dagger$ & $\mathbf{15}$ & 2 & 4 & 9 & $9_{45}^\dagger$ & $\mathbf{13}$ & 3 & 4 & 6\\
$6_{2}^\dagger$ & $\mathbf{6}$ & 2 & 2 & 2 & $9_{8}^\dagger$ & $\mathbf{10}$ & 2 & 3 & 5 & $9_{46}^\dagger$ & $\mathbf{14}$ & 2 & 5 & 7\\
$6_{3}$ & $\mathbf{6}$ & 2 & 2 & 2 & $9_{9}^\dagger$ & $\mathbf{22}$ & 3 & 4 & 15 & $9_{47}^\dagger$ & $\mathbf{9}$ & 2 & 3 & 4\\
$7_{1}^*$ & $\mathbf{17}$ & 3 & 5 & 9 & $9_{10}^\dagger$ & $\mathbf{21}$ & 3 & 3 & 15 & $9_{48}^\dagger$ & $\mathbf{10}$ & 3 & 3 & 4\\
$7_{2}^\dagger$ & $\mathbf{9}$ & 3 & 3 & 3 & $9_{11}^\dagger$ & $\mathbf{15}$ & 3 & 3 & 9 & $9_{49}^\dagger$ & $\mathbf{18}$ & 6 & 6 & 6\\
$7_{3}^*$ & $\mathbf{14}$ & 5 & 7 & 2 & $9_{12}^\dagger$ & $\mathbf{10}$ & 2 & 3 & 5 & $10_{1}^\dagger$ & $\mathbf{13}$ & 2 & 4 & 7\\
$7_{4}^\dagger$ & $\mathbf{11}$ & 2 & 2 & 7 & $9_{13}^\dagger$ & $\mathbf{19}$ & 3 & 3 & 13 & $10_{2}^\dagger$ & $\mathbf{20}$ & 4 & 4 & 12\\
$7_{5}^\dagger$ & $\mathbf{12}$ & 2 & 3 & 7 & $9_{14}^\dagger$ & $\mathbf{10}$ & 2 & 3 & 5 & $10_{3}^\dagger$ & $\mathbf{17}$ & 2 & 3 & 12\\
$7_{6}^\dagger$ & $\mathbf{7}$ & 2 & 2 & 3 & $9_{15}^\dagger$ & $\mathbf{10}$ & 2 & 3 & 5 & $10_{4}^\dagger$ & $\mathbf{16}$ & 2 & 6 & 8\\
$7_{7}$ & $\mathbf{7}$ & 2 & 3 & 2 & $9_{17}^\dagger$ & $\mathbf{9}$ & 3 & 3 & 3 & $10_{5}^\dagger$ & $\mathbf{20}$ & 4 & 4 & 12\\
$8_{1}^\dagger$ & $\mathbf{10}$ & 2 & 3 & 5 & $9_{18}^\dagger$ & $\mathbf{17}$ & 2 & 4 & 11 & $10_{6}^\dagger$ & $\mathbf{15}$ & 3 & 4 & 8\\
$8_{2}^\dagger$ & $\mathbf{12}$ & 3 & 3 & 6 & $9_{19}^\dagger$ & $\mathbf{10}$ & 2 & 3 & 5 & $10_{7}^\dagger$ & $\mathbf{12}$ & 2 & 4 & 6\\
$8_{3}^\dagger$ & $\mathbf{12}$ & 2 & 2 & 8 & $9_{20}^\dagger$ & $\mathbf{9}$ & 3 & 3 & 3 & $10_{8}^\dagger$ & $\mathbf{19}$ & 3 & 6 & 10\\
$8_{4}^\dagger$ & $\mathbf{11}$ & 2 & 4 & 5 & $9_{21}^\dagger$ & $\mathbf{12}$ & 2 & 3 & 7 & $10_{9}^\dagger$ & $\mathbf{18}$ & 5 & 5 & 8\\
$8_{5}^\dagger$ & $\mathbf{23}$ & 7 & 7 & 9 & $9_{22}^\dagger$ & $\mathbf{20}$ & 6 & 7 & 7 & $10_{10}$ & $\mathbf{12}$ & 4 & 6 & 2\\
$8_{6}^\dagger$ & $\mathbf{9}$ & 2 & 3 & 4 & $9_{23}^\dagger$ & $\mathbf{15}$ & 3 & 3 & 9 & $10_{11}^\dagger$ & $\mathbf{16}$ & 3 & 4 & 9\\
$8_{7}$ & $\mathbf{12}$ & 4 & 5 & 3 & $9_{24}^\dagger$ & $\mathbf{21}$ & 5 & 8 & 8 & $10_{12}^\dagger$ & $\mathbf{17}$ & 3 & 4 & 10\\
$8_{8}^\dagger$ & $\mathbf{9}$ & 2 & 3 & 4 & $9_{25}^\dagger$ & $\mathbf{22}$ & 5 & 7 & 10 & $10_{13}^\dagger$ & $\mathbf{15}$ & 2 & 3 & 10\\
$8_{9}^\dagger$ & $\mathbf{12}$ & 3 & 4 & 5 & $9_{26}^\dagger$ & $\mathbf{9}$ & 3 & 3 & 3 & $10_{14}^\dagger$ & $\mathbf{17}$ & 3 & 4 & 10\\
$8_{10}^\dagger$ & $\mathbf{21}$ & 6 & 7 & 8 & $9_{27}^\dagger$ & $\mathbf{9}$ & 3 & 3 & 3 & $10_{15}^\dagger$ & $\mathbf{15}$ & 3 & 4 & 8\\
$8_{11}^\dagger$ & $\mathbf{9}$ & 2 & 3 & 4 & $9_{28}^\dagger$ & $\mathbf{21}$ & 6 & 7 & 8 & $10_{16}^\dagger$ & $\mathbf{16}$ & 3 & 4 & 9\\
$8_{12}^\dagger$ & $\mathbf{10}$ & 2 & 2 & 6 & $9_{29}^\dagger$ & $\mathbf{22}$ & 5 & 8 & 9 & $10_{17}^\dagger$ & $\mathbf{18}$ & 5 & 5 & 8\\
$8_{13}$ & $\mathbf{9}$ & 2 & 4 & 3 & $9_{30}^\dagger$ & $\mathbf{20}$ & 5 & 6 & 9 & $10_{18}^\dagger$ & $\mathbf{14}$ & 2 & 6 & 6\\
$8_{14}^\dagger$ & $\mathbf{9}$ & 2 & 3 & 4 & $9_{31}$ & $\mathbf{9}$ & 3 & 3 & 3 & $10_{19}^\dagger$ & $\mathbf{17}$ & 3 & 4 & 10\\
$8_{15}$ & $\mathbf{27}$ & 10 & 5 & 12 & $9_{32}^\dagger$ & $\mathbf{11}$ & 3 & 3 & 5 & $10_{20}^\dagger$ & $\mathbf{12}$ & 2 & 4 & 6\\
$8_{16}^\dagger$ & $\mathbf{10}$ & 2 & 3 & 5 & $9_{33}^\dagger$ & $\mathbf{11}$ & 2 & 4 & 5 & $10_{21}^\dagger$ & $\mathbf{17}$ & 3 & 4 & 10\\
$8_{17}^\dagger$ & $\mathbf{10}$ & 2 & 4 & 4 & $9_{34}^\dagger$ & $\mathbf{10}$ & 2 & 4 & 4 & $10_{22}^\dagger$ & $\mathbf{17}$ & 3 & 6 & 8\\
$8_{18}^\dagger$ & $\mathbf{8}$ & 2 & 3 & 3 & $9_{35}^\dagger$ & $\mathbf{18}$ & 5 & 6 & 7 & $10_{23}^\dagger$ & $\mathbf{17}$ & 3 & 4 & 10\\
$8_{19}$ & $30$ & 5 & 7 & 18 & $9_{36}^\dagger$ & $\mathbf{20}$ & 6 & 7 & 7 & $10_{24}^\dagger$ & $\mathbf{12}$ & 2 & 4 & 6\\
$8_{20}^\dagger$ & $\mathbf{12}$ & 3 & 4 & 5 & $9_{37}^\dagger$ & $\mathbf{10}$ & 3 & 3 & 4 & $10_{25}^\dagger$ & $\mathbf{15}$ & 3 & 4 & 8\\
$8_{21}^\dagger$ & $\mathbf{12}$ & 3 & 4 & 5 & $9_{38}^\dagger$ & $\mathbf{24}$ & 6 & 7 & 11 & $10_{26}^\dagger$ & $\mathbf{17}$ & 3 & 6 & 8\\
$9_{1}^\dagger$ & $\mathbf{27}$ & 4 & 4 & 19 & $9_{39}^\dagger$ & $\mathbf{12}$ & 3 & 3 & 6 & $10_{27}^\dagger$ & $\mathbf{15}$ & 3 & 4 & 8\\
$9_{2}$ & $\mathbf{12}$ & 4 & 7 & 1 & $9_{40}^\dagger$ & $\mathbf{14}$ & 3 & 5 & 6 & $10_{28}^\dagger$ & $\mathbf{14}$ & 3 & 3 & 8\\
\newpage
$10_{29}^\dagger$ & $\mathbf{14}$ & 3 & 4 & 7 & $10_{77}^\dagger$ & $\mathbf{24}$ & 6 & 8 & 10 & $10_{121}^\dagger$ & $\mathbf{15}$ & 3 & 5 & 7\\
$10_{30}^\dagger$ & $\mathbf{14}$ & 3 & 3 & 8 & $10_{82}^\dagger$ & $\mathbf{15}$ & 4 & 4 & 7 & $10_{122}^\dagger$ & $\mathbf{10}$ & 3 & 3 & 4\\
$10_{31}^\dagger$ & $\mathbf{12}$ & 2 & 4 & 6 & $10_{83}^\dagger$ & $\mathbf{15}$ & 2 & 5 & 8 & $10_{123}^\dagger$ & $\mathbf{10}$ & 3 & 3 & 4\\
$10_{32}^\dagger$ & $\mathbf{13}$ & 4 & 4 & 5 & $10_{84}^\dagger$ & $\mathbf{13}$ & 3 & 3 & 7 & $10_{125}^\dagger$ & $\mathbf{26}$ & 8 & 9 & 9\\
$10_{33}^\dagger$ & $\mathbf{14}$ & 3 & 3 & 8 & $10_{85}^\dagger$ & $\mathbf{18}$ & 3 & 5 & 10 & $10_{129}^\dagger$ & $\mathbf{17}$ & 4 & 4 & 9\\
$10_{34}^\dagger$ & $\mathbf{12}$ & 2 & 4 & 6 & $10_{86}^\dagger$ & $\mathbf{13}$ & 2 & 5 & 6 & $10_{130}^\dagger$ & $\mathbf{19}$ & 3 & 8 & 8\\
$10_{35}^\dagger$ & $\mathbf{13}$ & 2 & 3 & 8 & $10_{87}^\dagger$ & $\mathbf{13}$ & 3 & 4 & 6 & $10_{131}^\dagger$ & $\mathbf{19}$ & 3 & 7 & 9\\
$10_{36}^\dagger$ & $\mathbf{12}$ & 2 & 4 & 6 & $10_{88}^\dagger$ & $\mathbf{12}$ & 2 & 5 & 5 & $10_{132}^\dagger$ & $\mathbf{15}$ & 3 & 5 & 7\\
$10_{37}^\dagger$ & $\mathbf{12}$ & 2 & 5 & 5 & $10_{89}^\dagger$ & $\mathbf{14}$ & 3 & 5 & 6 & $10_{133}^\dagger$ & $\mathbf{15}$ & 3 & 5 & 7\\
$10_{38}^\dagger$ & $\mathbf{12}$ & 2 & 5 & 5 & $10_{90}^\dagger$ & $\mathbf{16}$ & 3 & 4 & 9 & $10_{135}^\dagger$ & $\mathbf{23}$ & 4 & 9 & 10\\
$10_{39}^\dagger$ & $\mathbf{15}$ & 3 & 4 & 8 & $10_{91}^\dagger$ & $\mathbf{18}$ & 4 & 5 & 9 & $10_{136}^\dagger$ & $\mathbf{13}$ & 3 & 4 & 6\\
$10_{40}^\dagger$ & $\mathbf{15}$ & 3 & 4 & 8 & $10_{92}^\dagger$ & $\mathbf{19}$ & 6 & 9 & 4 & $10_{137}^\dagger$ & $\mathbf{14}$ & 4 & 4 & 6\\
$10_{41}^\dagger$ & $\mathbf{10}$ & 3 & 3 & 4 & $10_{93}^\dagger$ & $\mathbf{17}$ & 5 & 6 & 6 & $10_{138}^\dagger$ & $\mathbf{20}$ & 6 & 7 & 7\\
$10_{42}^\dagger$ & $\mathbf{10}$ & 3 & 3 & 4 & $10_{94}^\dagger$ & $\mathbf{18}$ & 5 & 6 & 7 & $10_{140}^\dagger$ & $\mathbf{21}$ & 6 & 7 & 8\\
$10_{43}^\dagger$ & $\mathbf{10}$ & 3 & 3 & 4 & $10_{95}^\dagger$ & $\mathbf{16}$ & 4 & 5 & 7 & $10_{141}^\dagger$ & $\mathbf{14}$ & 4 & 5 & 5\\
$10_{44}^\dagger$ & $\mathbf{10}$ & 3 & 3 & 4 & $10_{97}^\dagger$ & $\mathbf{20}$ & 3 & 5 & 12 & $10_{143}^\dagger$ & $\mathbf{20}$ & 4 & 5 & 11\\
$10_{45}^\dagger$ & $\mathbf{10}$ & 3 & 3 & 4 & $10_{102}^\dagger$ & $\mathbf{15}$ & 4 & 4 & 7 & $10_{144}^\dagger$ & $\mathbf{17}$ & 4 & 6 & 7\\
$10_{52}^\dagger$ & $\mathbf{24}$ & 6 & 8 & 10 & $10_{103}^\dagger$ & $\mathbf{15}$ & 3 & 3 & 9 & $10_{145}^\dagger$ & $\mathbf{14}$ & 3 & 5 & 6\\
$10_{57}^\dagger$ & $\mathbf{26}$ & 8 & 9 & 9 & $10_{104}^\dagger$ & $\mathbf{14}$ & 3 & 5 & 6 & $10_{146}^\dagger$ & $\mathbf{10}$ & 2 & 3 & 5\\
$10_{58}^\dagger$ & $\mathbf{23}$ & 7 & 7 & 9 & $10_{105}^\dagger$ & $\mathbf{15}$ & 4 & 4 & 7 & $10_{147}^\dagger$ & $\mathbf{10}$ & 3 & 3 & 4\\
$10_{59}^\dagger$ & $\mathbf{19}$ & 5 & 7 & 7 & $10_{106}^\dagger$ & $\mathbf{17}$ & 4 & 5 & 8 & $10_{148}^\dagger$ & $\mathbf{18}$ & 4 & 5 & 9\\
$10_{60}^\dagger$ & $\mathbf{20}$ & 6 & 7 & 7 & $10_{107}^\dagger$ & $\mathbf{12}$ & 3 & 4 & 5 & $10_{149}^\dagger$ & $\mathbf{18}$ & 4 & 5 & 9\\
$10_{65}^\dagger$ & $\mathbf{24}$ & 6 & 6 & 12 & $10_{108}^\dagger$ & $\mathbf{15}$ & 3 & 3 & 9 & $10_{150}^\dagger$ & $\mathbf{15}$ & 3 & 6 & 6\\
$10_{67}^\dagger$ & $\mathbf{16}$ & 2 & 6 & 8 & $10_{109}^\dagger$ & $\mathbf{26}$ & 6 & 8 & 12 & $10_{151}^\dagger$ & $\mathbf{11}$ & 3 & 3 & 5\\
$10_{68}^\dagger$ & $\mathbf{14}$ & 4 & 4 & 6 & $10_{110}^\dagger$ & $\mathbf{16}$ & 4 & 5 & 7 & $10_{155}^\dagger$ & $\mathbf{18}$ & 4 & 7 & 7\\
$10_{69}^\dagger$ & $\mathbf{14}$ & 3 & 5 & 6 & $10_{112}^\dagger$ & $\mathbf{10}$ & 3 & 3 & 4 & $10_{156}^\dagger$ & $\mathbf{10}$ & 3 & 3 & 4\\
$10_{70}^\dagger$ & $\mathbf{24}$ & 5 & 8 & 11 & $10_{113}^\dagger$ & $\mathbf{10}$ & 3 & 3 & 4 & $10_{157}^\dagger$ & $\mathbf{14}$ & 2 & 4 & 8\\
$10_{71}^\dagger$ & $\mathbf{20}$ & 6 & 7 & 7 & $10_{114}^\dagger$ & $\mathbf{11}$ & 2 & 4 & 5 & $10_{158}^\dagger$ & $\mathbf{12}$ & 2 & 5 & 5\\
$10_{72}^\dagger$ & $\mathbf{24}$ & 6 & 8 & 10 & $10_{115}^\dagger$ & $\mathbf{14}$ & 2 & 6 & 6 & $10_{159}^\dagger$ & $\mathbf{11}$ & 3 & 4 & 4\\
$10_{73}^\dagger$ & $\mathbf{20}$ & 6 & 7 & 7 & $10_{116}^\dagger$ & $\mathbf{16}$ & 4 & 5 & 7 & $10_{160}^\dagger$ & $\mathbf{11}$ & 2 & 3 & 6\\
$10_{74}^\dagger$ & $\mathbf{14}$ & 4 & 4 & 6 & $10_{117}^\dagger$ & $\mathbf{15}$ & 3 & 6 & 6 & $10_{161}^\dagger$ & $\mathbf{23}$ & 6 & 7 & 10\\
$10_{75}^\dagger$ & $\mathbf{14}$ & 3 & 5 & 6 & $10_{118}^\dagger$ & $\mathbf{12}$ & 3 & 4 & 5 & $10_{163}^\dagger$ & $\mathbf{11}$ & 2 & 4 & 5\\
$10_{76}^\dagger$ & $\mathbf{26}$ & 5 & 8 & 13 & $10_{119}^\dagger$ & $\mathbf{13}$ & 3 & 5 & 5 & $10_{164}^\dagger$ & $\mathbf{14}$ & 4 & 4 & 6\\
\end{longtable}}

\begin{figure}[p]
\centering
\includegraphics[height=0.84\textheight]{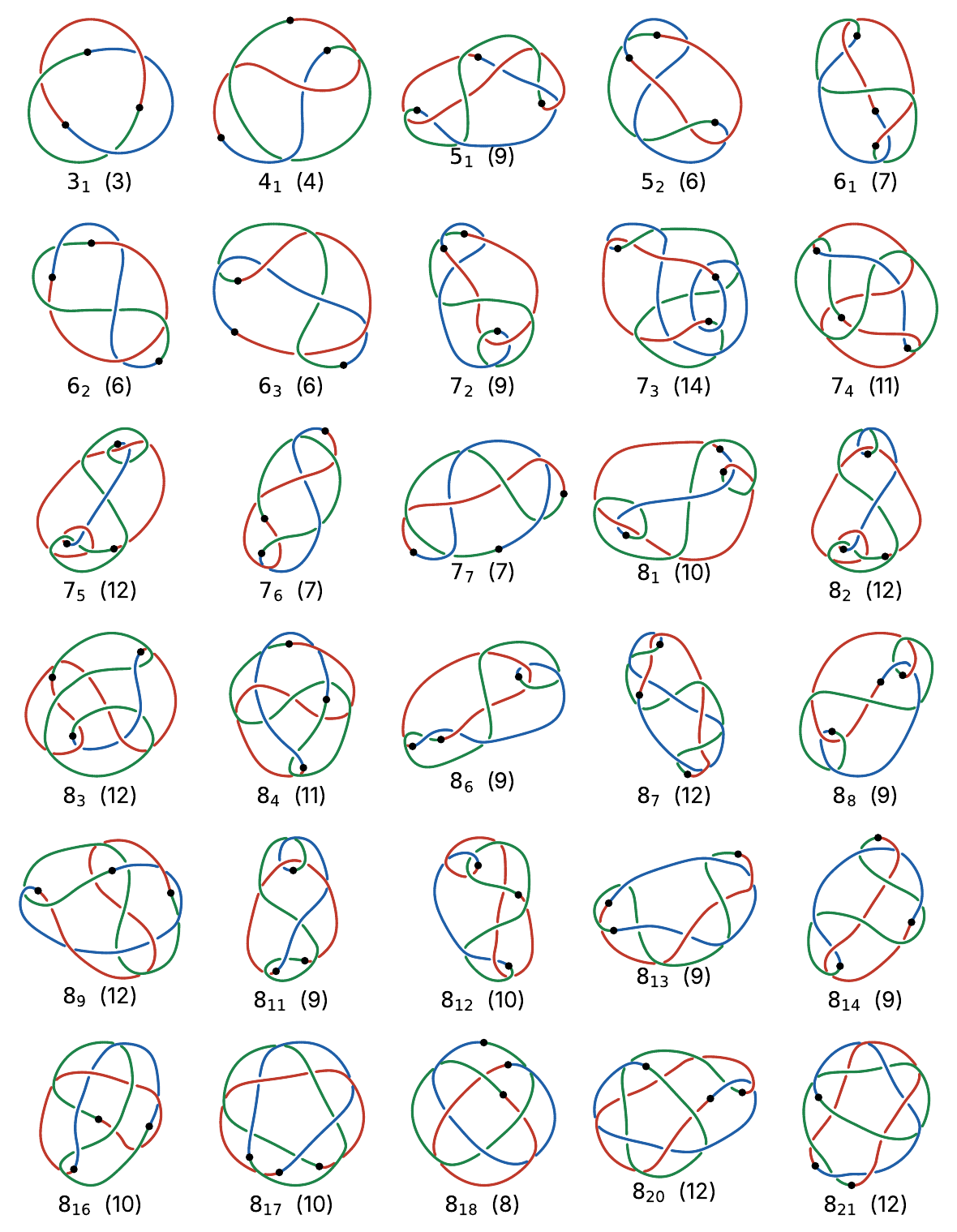}
\caption{The $3$-segment presentations of Table~\ref{tab:census} of the $30$ prime knots with at most eight crossings other than $7_1$, $8_5$, $8_{10}$, $8_{15}$ and $8_{19}$ (Figure~\ref{fig:gallery-large}), with the number of crossings in parentheses. The segments $e_1$ (red), $e_2$ (blue) and $e_3$ (green) are separated by the three black dots, and $e_1>e_2>e_3>e_1$ at every crossing. Each drawing has the signed Gauss code of the presentation in the supplementary data.}
\label{fig:gallery}
\end{figure}

\begin{figure}[p]
\centering
\includegraphics[width=\textwidth]{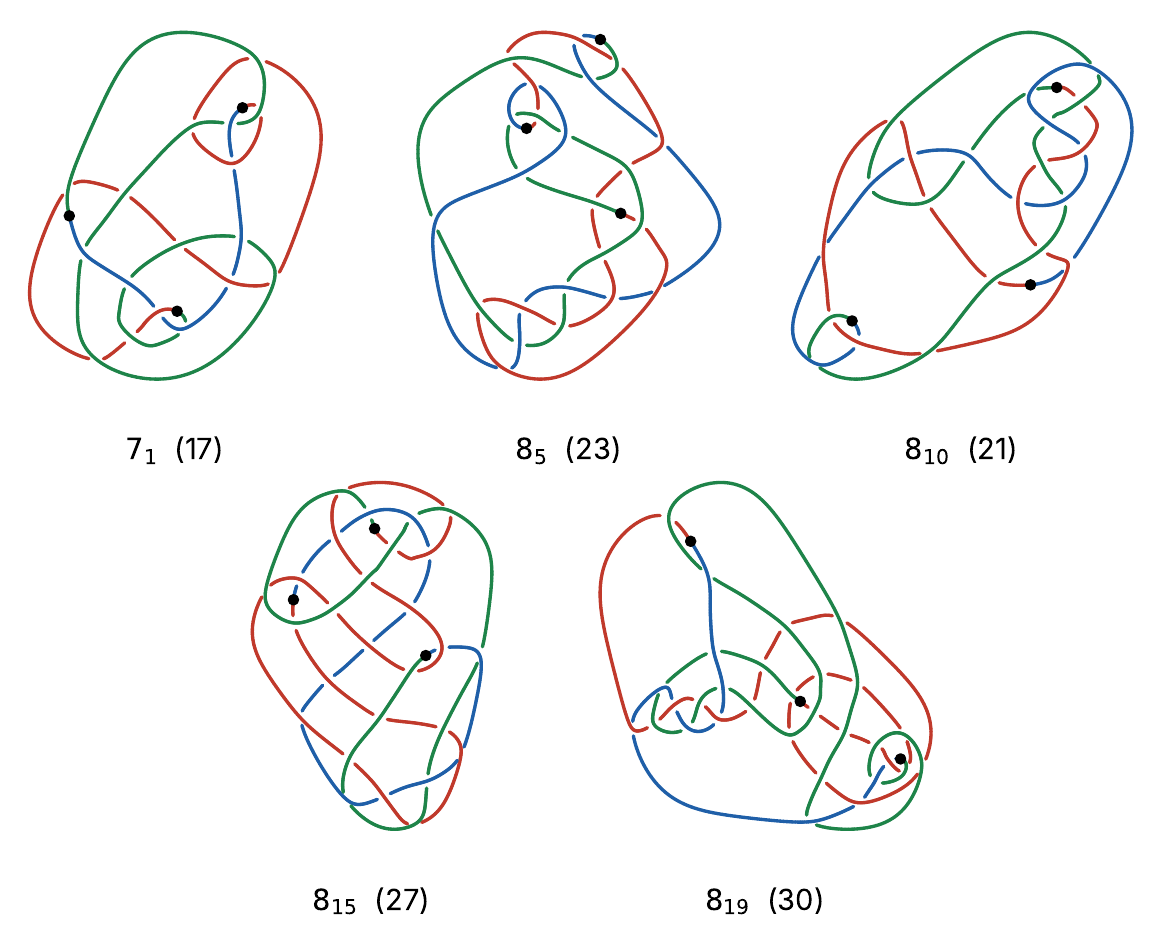}
\caption{The $3$-segment presentations of Table~\ref{tab:census} of the five prime knots with at most eight crossings which need more than $14$ crossings: $7_1$, $8_5$, $8_{10}$, $8_{15}$ and $8_{19}$. Conventions as in Figure~\ref{fig:gallery}.}
\label{fig:gallery-large}
\end{figure}

\end{document}